\documentclass[12pt,reqno]{amsart}

\usepackage{amsmath}
\usepackage{amsthm}
\usepackage{amscd}
\usepackage{amssymb}
\usepackage{graphicx}
\usepackage{epsfig}
\theoremstyle{plain}
\newtheorem{theorem}{Theorem}[section]
\newtheorem{proposition}[theorem]{Proposition}
\newtheorem{lemma}[theorem]{Lemma}
\newtheorem{corollary}[theorem]{Corollary}

\newtheorem{definition}[theorem]{Definition}
\newtheorem{remark}[theorem]{Remark}

\def\int{\text{interior}}
\def\hyp {\hbox {\rm {H \kern -2.8ex I}\kern 1.25ex}}
\def\reals {\hbox {\rm {R \kern -2.8ex I}\kern 1.15ex}}
\def\integers {\hbox {\rm { Z \kern -2.8ex Z}\kern 1.15ex}}
\def\naturals {\hbox {\rm {N \kern -2.8ex I}\kern 1.20ex}}
\def\rationals {\hbox {\rm { Q \kern -2.2ex l}\kern 1.15ex}}
\def\hyp {\hbox {\rm {H \kern -2.7ex I}\kern 1.25ex}}
\def\bar{\overline}

\def\hat{\widehat}

\def\strutdepth{\dp\strutbox}
\def \ss{\strut\vadjust{\kern-\strutdepth \sss}}
\def \sss{\vtop to \strutdepth{
\baselineskip\strutdepth\vss\llap{$\diamondsuit\;\;$}\null}}

\def\strutdepth{\dp\strutbox}
\def \sst{\strut\vadjust{\kern-\strutdepth \ssss}}
\def \ssss{\vtop to \strutdepth{
\baselineskip\strutdepth\vss\llap{$\spadesuit\;\;$}\null}}

\def\strutdepth{\dp\strutbox}
\def \ssh{\strut\vadjust{\kern-\strutdepth \sssh}}
\def \sssh{\vtop to \strutdepth{
\baselineskip\strutdepth\vss\llap{$\heartsuit\;\;$}\null}}

\begin{document}

 \title[Horizontal Dehn surgery and genericity]{Horizontal Dehn surgery and genericity in the curve complex}

\author{Martin Lustig} 
 
\author{Yoav Moriah}

\thanks {This research is supported by  a grant from the High Council 
for Scientific and Technological Cooperation between France and 
Israel}

\begin{abstract}    

We introduce a general notion of ``genericity" for countable subsets of a space with Borel measure, 
and apply it to the set of  vertices  in the curve  complex  of a surface $\Sigma$,  interpreted as subset 
of the space of projective measured laminations in $\Sigma$, equipped with its natural Lebesgue measure. 

We prove that,  for any $3$-manifold $M$, the set of curves $c$ on a Heegaard surface 
$\Sigma \subset M$, such that every non-trivial Dehn twist at $c$ yields a  Heegaard splitting of high distance, is generic in the set of all essential simple closed curves on $\Sigma$. 

Our definition of ``genericity" is different  and more intrinsic than alternative such existing notions, 
given e.g. via random walks or via limits of quotients of finite sets.

\end{abstract}


\address{Math\'ematiques (LATP)\\
Universit\'e P. C\'ezanne - Aix-Marseille III \\
Ave. Escad. Normandie-Niemen, 13397 Marseille 20, France}
\email{martin.lustig@univ-cezanne.fr}

\address{Department of Mathematics \\
Technion \\
Haifa, 32000 Israel}
\email{ymoriah@tx.technion.ac.il}

\keywords{Train tracks, Curve complex, Heegaard distance, Heegaard 
splittings}

\maketitle

\section{Introduction} \label{introduction}

\vskip10pt

For any compact connected orientable surface $\Sigma$ of genus $g \geq 2$ the {\em curve complex},  
denoted by $\mathcal{C}(\Sigma)$, is a locally infinite  simplicial complex, where every $m$-simplex is  determined  by a collection of $m+1$ isotopy classes of pairwise disjoint simple closed essential curves on  $\Sigma$. Every handlebody $W$ with  boundary identification map 
$h: \partial W \overset{\approx}{\longrightarrow} \Sigma$ determines a  {\it  disk complex} $\mathcal{D}(W) \subset \mathcal{C}(\Sigma)$, which is defined by  the condition that every curve  representing a 
vertex of $\mathcal{D}(W)$  must bound a disk in $W$. 
 
For any simple closed essential curve $c \subset \Sigma$ we denote by $\delta_c$ the Dehn twist at $c$.
One can perturb the identification map $h: \partial W \overset{\approx}{\longrightarrow}  \Sigma$
by composing it with a power $\delta_c^m$, to get a  handlebody $W_c^m$ which is a homeomorphic copy of $W$ (via some $g:  W_c^m \overset{\approx}{\longrightarrow} W$) and has 
$\delta_c^m\circ h \circ g: \partial W_c^m  \to \Sigma$  as the induced boundary identification map.

Any Heegaard splitting $M = V \underset{\partial V \approx \Sigma \approx \partial W}{\cup} W$ 
of a compact closed orientable 3-manifold $M$ defines two subcomplexes  
$\mathcal{D}(V), \mathcal{D}(W) \subset \mathcal{C}(\Sigma)$, and the minimal distance (in the simplicial metric of $\mathcal{C}(\Sigma)$)  between points in these subcomplexes is called the 
{\it distance} of the Heegaard splitting, denoted $d(V, W)$.  

Thus any simple closed essential curve  $c \subset \Sigma \subset M$ and any integer 
$m \in \mathbb Z$ define  a new 3-manifold $M_c^m$ with Heegaard splitting: 
$$M_c^m = V \underset{\partial V \approx \Sigma \approx \partial W_c^m}{\cup} W_c^m$$
The curve $c$ is called {\it $n$-optimal} if for  every non-trivial  twist exponent $m \neq 0$ one has: 
$$d(V, W^m_c) \geq n$$

\vskip5pt

The main result of this paper,  Theorem \ref{twistgenericity}, can be paraphrased as follows:

\vskip5pt 

\begin{theorem} \label{optimal}
Let  $M$ be a closed orientable $3$-manifold, and let $\Sigma$ be a Heegaard surface of $M$.  Then 
almost every essential simple closed curve $c$ on  $\Sigma$  is  $n$-optimal, for any $n \geq 1$.

\end{theorem}

\vskip5pt

Here the terminology ``almost every'' refers to the Lebesgue measure on $\mathcal{PML}(\Sigma)$, 
the space of projective measured laminations on $\Sigma$.  Typically, genericity results come from  random walks and almost always involve at the very beginning the  choice of extra data, for example a preferred generating system of a group (see e.g.~\cite{Du},  ~\cite{Mah}).   Generic sets tend to vary if one modifies these extra data.

Alternatively, a ``complexity'' function with finite preimage subsets is used, and genericity is defined by stating that the quotient of the cardinalities of certain sets of bounded complexity tends to 1 if the complexity bound tends to $\infty$  (see e.g. ~\cite{Ar}, ~\cite{Ol}, ~\cite{KS}). Again, such genericity results depend heavily on the choice of the complexity function at the very beginning.
 
The concept of {\it genericity} introduced in this paper is independent of any such additional choices and hence is somewhat preferable.  We give a detailed discussion of this concept in Section \ref{generic} of this paper.

\vskip7pt

\begin{remark}\label{threesphere}\rm

The proof of the above theorem is constructive. For example, in the special case of the 3-sphere 
$M = S^3$ with a standard Heegaard surface $\Sigma$ of genus $g$, it exhibits, for any  given 
$ n \geq 1$,  explicit curves $c$ on $\Sigma$ with the property that a single Dehn twist at $c$ alters the standard Heegaard splitting of $S^3$ to a new  Heegaard splittings (of the new manifold  $S_c^{3,1}$) with  distance bigger or equal  to $n$. Notice also that the manifolds $S_c^{3,m}$ can alternatively be obtained by $\frac{1}{m}$-surgery on  the  knot  $c \subset S^{3}$,  where the surgery coefficients are defined with respect to the meridian on  $\partial N(c)$ and the ``horizontal'' slope $\Sigma  \cap  \partial N(c)$. Manifolds $M_c^m$ defined as above are said to be obtained from $M$ by {\it horizontal Dehn surgery on the knot $c \subset \Sigma$} (see e.g. \cite{MS}).

\end{remark}

\vskip7pt

The curve complex associated with a closed surface has become over the past ten years a subject 
of increasing importance for low dimensional topology.  It is known to be a $\delta$-hyperbolic 
space in the sense of Gromov (see ~\cite{MM}).  For a $3$-manifold $M$ with Heegaard splitting 
$M = V \cup_{\Sigma} W$ the above defined distance $d(V, W)$ (sometimes also referred to as the 
{\it Hempel distance} of the splitting)  has become an important invariant of the Heegaard splitting.  
For example, it has been shown  that any  $3$-manifold  admits only finitely many Heegaard splittings
of distance $\geq 3$  (see ~\cite{Sc1}).  Furthermore, it follows from Perelman's proof of the Geometrization  Conjecture and from the classification of Heegaard splittings of Seifert fibered spaces (see ~\cite{MS} and ~\cite{He}) that every 3-manifold $M$ with at least one Heegaard splitting of distance $\geq 3$ is hyperbolic.

In the process of proving the above results, we have also derived the following two genericity statements about distance in the curve complex, which may be of interest in their own right.  Since they confirm what most experts feel ought to be true, they can alternatively be viewed as confirmation that the definition of ``genericity'' introduced in this paper is a useful and natural notion.  We have:

\vskip10pt

\noindent{\bf Corollary  \ref {longdistancefromc}.}
For any essential simple closed curve $c$ on $\Sigma$, the set $\mathcal C^{0}_{n}(c)$ of all essential simple closed curves $k$ on $\Sigma$ with distance $$d(k, c) \geq n$$ is generic in the set 
$\mathcal C^{0}(\Sigma)$ of all essential simple closed curves on $\Sigma$.

\vskip15pt

\noindent{\bf Theorem \ref{longdistancefromH}. }
For any handlebody $H$ with boundary surface $\partial H = \Sigma$ the set $\mathcal C^{0}_{n}(H)$ 
of all essential simple closed curves $k$ on $\Sigma$ with distance
$$d(k, H)  = \min\{d(k, c) \mid c \in \mathcal{D}(H)\} \geq n$$ 
is generic in the set  $\mathcal  C^{0}(\Sigma)$.

\vskip10pt

We would like to point out that an important ingredient in the proofs of the above theorems is 
Kerckhoff's result that the limit set of the  handlebody group has measure zero in the Thurston 
boundary (see ~\cite{Ke}).

\vskip15pt

\noindent {\bf Acknowledgments.} 
We thank  C. Connoll,  I. Kapovitch,  Y. Minsky\footnote{ \/ Y. Minsky  also provided considerable help with the figures.}   and P. Hubert for helpful remarks. 
Furthermore we  would like to thank the High Council for Scientific and Technological Cooperation  between France  and Israel for its support, also the Universit\'e de Provence in  Marseilles and the Technion in Haifa. The second author would especially like to thank Yale University for its hospitality.

\vskip40pt

\section{Notation and background}\label{sectionone}

\vskip4pt

In this section we will recall various definitions and background  material needed for the  following sections. Most of the material of this section has been presented in full detail  in  ~\cite {LM1}, and is here only recalled briefly, for the convenience of the reader. 

\vskip10pt

\subsection{The curve complex}\label{The curve complex}\hfill

\vskip10pt

Given an orientable connected surface $\Sigma$ of genus  $g \geq 2$,  the {\it curve  complex},  denoted by $\mathcal {C}(\Sigma)$  is defined as follows:

\begin{enumerate}
\item The set of vertices $\mathcal{C}^0(\Sigma)$ is the set of   isotopy classes of simple closed 
curves on $\Sigma$.

\item  An $n$-simplex in $\mathcal{C}(\Sigma)$ is a collection $\{v_0, ..., v_{n}\}$ of    vertices which 
can be represented by mutually disjoint curves.
\end{enumerate}

On the $1$-skeleton of $\mathcal {C}^1 (\Sigma)$  one defines a metric $d_{\mathcal{C}}( \cdot,\cdot)$  
by declaring the length of every edge to be $1$. For the purpose of  this paper it will suffice to consider 
only  $\mathcal {C}^1 (\Sigma)$. 

\vskip10pt 

\subsection{Train tracks}\hfill

\vskip10pt

A {\it train track} $\tau$ in $\Sigma$ is a compact  subsurface with boundary, which is equipped with a  {\it singular $I$-fiberation}: The interior of  $\tau$ is fibered by open arcs,   and   the  fibration extends to a fiberation of the compact surface  $\tau$  by  properly  embedded closed arcs (the {\it $I$-fibers}), except for finitely many {\it singular points}  (also called {\it cusp points}) on $\partial \tau$,  where precisely two  fibers meet. We call these  fibers {\it singular fibers}. We admit  the case that a fiber is {\it  doubly singular}, i.e. both of  its  endpoints are singular points. 

Two singular fibers are {\it adjacent} if they share a singular  point  as a common endpoint. A maximal connected union of singular  or  doubly singular $I$-fibers is called an {\it exceptional fiber}.  It is either 
homeomorphic  to a closed  interval, or to a simple closed curve on  $\Sigma$.  In  the latter case it will be called a {\it  cyclic  exceptional fiber}. We explicitely admit this second case, although  in the classical  train track  literature this case is sometimes  suppressed.

\vskip5pt

\begin{definition}\label{fatraintrack}\rm 
A train track $\tau \subset \Sigma$ is called {\it fat} if all of  its  exceptional fibers are cyclic. We denote by  ${\mathcal E}_{\tau}$ the collection of simple closed curves on    $\Sigma$ given by the exceptional fibers of  $\tau$.
\end{definition}

A train track $\tau$ in $\Sigma$ is called {\it filling}, if all  complementary components of $\tau$ in 
$\Sigma $ are simply  connected. The train track $\tau$ is called {\it maximal}, if every  complementary component is a {\it triangle}, i.e. it is simply  connected and there are precisely three singular points on its boundary.

\vskip8pt

An arc, a closed curve or a lamination in $\Sigma$ {\it is carried}  by a train  track  $\tau \subset \Sigma$ if it is contained in $\tau$ and  is throughout transverse to the  $I$-fibers of $\tau$.   Two simple arcs 
carried by $\tau$ are {\it parallel} if  they intersect  the same $I$-fibers,  and these intersections occur on the two arcs  in precisely the same order. An arc, a closed curve or a lamination on $\Sigma$ which is  carried by $\tau$ is   said to  {\it cover} $\tau$ if  it  meets every $I$-fiber of $\tau$.

\vskip10pt

\subsection{Unzipping paths and derived train tracks}\hfill
\label{unzippingplusderived}

\vskip10pt

Given a train track $\tau \subset \Sigma$ which carries a  lamination $\mathcal{L}$ we can obtain 
a new train track, which still carries $\mathcal{L}$, as follows:

\vskip8pt

The train track $\tau$ can be {\it split} by moving any of the singular points $P$ (now called a 
{\it zipper}),  which is located on the boundary of a complementary component   $\Delta$ of $\tau$,   into the interior of $\tau$.   The   zipper $P$ will move along an {\it unzipping path}, which is embedded in  the  interior of  $\tau \smallsetminus \mathcal{L}$ and is transverse to the   $I$-fibers. Two unzipping paths are not allowed to cross each   other.  An unzipping path which covers  $\tau$ is called  
{\it complete}. 

In case two  zippers meet the same connected component of an  $I$-fiber in  
$\tau\smallsetminus {\mathcal L}$ from different directions, they  have to join up, thus changing the topology of the train track and of  its  complementary components. A situation like this is called a 
{\it  collision}.  In case of a collision the unzipping  procedure  stops. 

\vskip8pt

\begin{definition} \label{derivedtt} \rm 
We say that $\tau$ can be  derived with respect  to $\mathcal L$  if we can successively (or  simultaneously, it does not make any  difference) unzip every zipper along  a complete unzipping path, without ever running into a  collision.   The train  track $\tau'$  obtained by unzipping along  shortest possible  complete unzipping paths is  said to be  {\it  derived from $\tau$ with respect to $\mathcal L$}, or simply {\it  derived from  $\tau$}.  

\end{definition}

\vskip8pt

\begin{remark}\label{sametype} \rm 
If the train track $\tau'$ is derived from a maximal train track $\tau$, then   every complementary component $\Delta$ of $\tau'$ is also a triangle, i.e. $\tau'$ is also maximal. This follows directly from  Definition \ref{derivedtt},   since during the deriving process the unzipping paths never run into 
collisions. 

\end{remark}

\vskip8pt

\begin{lemma} [\cite{LM}, Lemma 2.9] \label{carriesimpliesunzip} 

Given a surface $\Sigma$ and  maximal train tracks  $\tau, \tau' \subset  \Sigma$  so that  $\tau'$ is derived from $\tau$. Let  $D$ be a simple closed  curve carried by $\tau'$. 
Then $D$  covers $\tau$.

\end{lemma}

\vskip8pt

A collection of train tracks  $\tau_{0} \supset \tau_{1} \supset   \ldots \supset\tau_{n}$ will be  
called an  {\it $n$-tower  of derived train tracks in $\Sigma$} if    each $\tau_{i}$ is  derived from  
$\tau_{i-1}$, for all $i  = 1, \ldots, n$. In this case we say that $\tau_{n}$ has been {\it $n$ times derived} 
from $\tau_{0}$.

\vskip12pt

\subsection{Complete fat train tracks} \label{fattraintracks} \hfill

\vskip12pt
 
A curve $D$  is called {\it tight} with   respect to a system of pairwise disjoint essential simple closed
curves $\mathcal E = \{E_{1}, \ldots, E_{r}\}$  in $\Sigma$  if the number of  intersection points with 
$\mathcal E$ can not be strictly decreased  by an isotopy of  $D$. The same terminology is used 
for arcs $\alpha$ which have their  endpoints on $\mathcal E$, where the endpoints cannot leave 
$\mathcal E$ throughout  the isotopy.

\vskip10pt

\begin{definition}\label{smallwave}\rm  
Let $P \subset \Sigma$ be a pair-of-pants,  i.e. a sphere with three open disks removed. 

\vskip7pt

\begin{enumerate}
\item[(a)]  A simple arc in $P$ which has its two endpoints on different  components of $\partial P$ 
will  be called a {\it seam}.

\item[(b)]  A simple arc in $P$ which has  both endpoints on the same  component of $\partial P$, 
and is not  $\partial$-parallel, will be called a {\it wave}.

\item[(c)]   An essential simple closed  curve $D \subset \Sigma$ has a  wave   (or a seam) 
with respect to a  system of curves $\mathcal{E}  \subset \Sigma$ if $D$ is tight with respect to 
$\mathcal E$ and if it contains a subarc  that is a wave (or a seam) in a  complementary component 
$P_{i}$ of  $\mathcal{E}$ in $\Sigma$ which is a pair-of-pants.  

\item[(d)]   An essential simple closed  curve $D \subset \Sigma$ has a {\it  wave} with respect to a fat  train  track  $\tau$ if $D$ has a wave with respect to ${\mathcal  E}_{\tau}$, or if $D$ is isotopic to some  $E_{k} \in {\mathcal E}_{\tau}$.
\end{enumerate}
\end{definition}

\vskip10pt

A system $\mathcal E$ of pairwise disjoint essential simple closed curves on $\Sigma$  is called  a 
{\it complete decomposing system} if every complementary  component of  $\mathcal E$ in $\Sigma$ is 
a  pair-of-pants.

\vskip10pt

\begin{definition} \label{completefatraintrack}\rm 
A fat train track $\tau \subset  \Sigma$ is called {\it complete}  if  the following conditions are satisfied:
\begin{enumerate}

\item The collection $\mathcal{E}_{\tau}$ of exceptional fibers of  $\tau$ is a complete decomposing system on  $\Sigma$.

\item Each pair-of-pants $P_{i}$ complementary to the system  $\mathcal{E}_{\tau}$ contains two triangles as complementary  components of $\tau$ in $P_{i}$. 

\item  The train track $\tau$ only carries seams, but no waves, with respect  to the complete decomposing system $\mathcal{E}_{\tau}$.

\end{enumerate}

\end{definition}

\vskip8pt

Notice that every complete fat train track is in particular maximal.

\vskip8pt

\begin{remark} \label{fillingpants} \rm
Let $\mathcal E$ be a complete decomposing system on the surface   $\Sigma$, and let $D$ be an essential simple closed curve (or a system of  such curves) on $\Sigma$ that is tight with respect to 
$\mathcal E$. We say that $D$ {\em fills}  a pair-of-pants $P$  complementary to $\mathcal E$, if 
$D \cap P$ is the disjoint union of precisely  3  distinct isotopy  classes of intersection arcs. 
Then the  following three statements are equivalent:

\begin{enumerate}
\item The curve $D$ fills every pair-of-pants complementary to  $\mathcal E$, and none of the intersection arcs is a wave.

\item There exists a unique complete fat train track $\tau$ with  exceptional fibers  
$\mathcal{E}_{\tau} = \mathcal{E}$ that carries $D$.

\item There exists some complete fat train track $\tau$ with   exceptional fibers 
$\mathcal{E}_{\tau} = \mathcal{E}$ that is covered  by $D$.

\end{enumerate}

\end{remark}

\vskip4pt

\begin{lemma}[Lemma 3.5 of \cite{LM1}]  \label{fatcarries}
Let  $\mathcal{E} \subset \Sigma$ be a  complete decomposing system.  Any essential simple closed curve $D \subset \Sigma$  which does not  have  a wave with respect to $\mathcal{E}$, and is not parallel to  any $E_{k} \in \mathcal{E}$, is carried by some complete fat train  track $\tau$ with exceptional fibers  ${\mathcal E}_{\tau} = {\mathcal E}$.

The same is true for any system  $\mathcal D$ of pairwise disjoint essential  simple  closed curves which satisfy the same conditions as $D$,  or any tight lamination $\mathcal L$.

\end{lemma}

\vskip10pt

\subsection{Tight train tracks} \label{tighttraintracks} \hfill

\vskip10pt

\begin{definition} \label{tighttts} \rm

A maximal train track $\tau$ is called {\it tight} with respect to some complete decomposing system 
$\mathcal E$ on $\Sigma$ if the following conditions are satisfied:

\begin{enumerate}

\item[(1)] For any curve $E_k \in \mathcal{E}$ every connected component of $E_k \cap \tau$ is a 
disjoint union of (possibly exceptional) $I$-fibers of $\tau$.

\item[(2)] For every connected component $\Delta_j$ complementary to $\tau$ the intersection segments with any $E_k  \in \mathcal{E}$ are arcs with endpoints on two distinct  sides of $\Delta_j$.

\item[(3)] Each of the three cusps of any complementary component $\Delta_j$ is contained in some of the $E_k$.

\end{enumerate}

\end{definition}

The condition (3) of Definition \ref{tighttts} is equivalent to stating that every singular $I$-fiber of $\tau$ lies on some of the curves $E_i \in \mathcal E$.

The reader may want to note that the above definition is less restrictive than what it appears:  A train track which is tight with respect to $\mathcal E$ may well carry a wave with respect to $\mathcal E$ !

\vskip10pt

\begin{lemma}
\label{fatimpliestight}
If a maximal train track $\tau$ is obtained from 
deriving finitely many times some fat train track $\hat\tau$ on $\Sigma$ with 
$\mathcal{E}_{\hat\tau} = \mathcal E$, then $\tau$ is tight with respect to $\mathcal E$. 

\end{lemma}

\vskip7pt

\begin{proof}
It follows directly from the definition of a fat train track (see Definition \ref{fattraintracks}) that $\hat \tau$ satisfies all three conditions of Definition \ref{tighttts}. On the other hand, one verifies directly that conditions (1) and (2) are preserved in the unzipping process. 

In order to prove condition (3), we assume by induction that it is satisfied by the train track $\tau'$ 
from which $\tau$ is derived. By Definition \ref{derivedtt} for each cusp point $Z$ of a complementary component $\Delta_j$ of $\tau$ there is a cusp point  $Z'$ on the complementary component 
$\Delta'_j \subset \Delta_j$ of $\tau'$ and an unzipping path $\sigma_i$ that starts at $Z'$ and ends 
in $Z$.  Recall that by definition of the deriving process no unzipping path can contain a proper initial subpath that meets every $I$-fiber of $\tau'$.

Now, since $\sigma_i$ is transverse to the $I$-fibering of $\tau'$, every time that $\sigma_i$ traverses 
an $I$-fiber $I_m$ of $\tau'$ contained in some $E_k$, there is a well defined adjacent (possibly exceptional) $I$-fiber $I_{m+1}$ contained in some $E_l$ which must be traversed next by $\sigma_i$: This fact  follows from our inductive hypothesis, which implies (via property (3) of Definition \ref{tighttts}) that every singular fiber of $\tau'$ is contained in some $E_h$.

As a consequence, if some $I$-fiber $I'$ of $\tau'$ between $I_m$ and $I_{m+1}$ has  been traversed by a proper subpath of $\sigma_i$, then the same must be true for any other $I$-fiber $I''$ between 
$I_m$ and $I_{m+1}$.  It follows that the endpoint $Z$ of $\sigma_i$ must be a singular fiber of $\tau'$, which  has been shown above to lie on some $E_k$.

\end{proof}

\begin{proposition}
\label{carriedcurves}
Let $\mathcal E$ be a complete decomposing system of $\Sigma$, and let $\tau$ be a maximal train track that is tight with respect to $\mathcal E$. Let $c$ be a simple closed curve  (or a finite collection of such) on $\Sigma$ that is tight with respect to $\mathcal E$ and contains a subarc $\beta$ which covers $\tau$.  Then $c$ can be carried by $\tau$.

\end{proposition}

\begin{proof}
We first consider the cyclic sequence of intersection points of $c$ with the singular fibers of $\tau$, and the segments $c_i$ between any two of them, to get a decomposition of $c$ as cyclic concatenation of the $c_i$. If any of the $c_i$ is contained in $\tau$, we may assume after an inessential isotopy that $c_i$ is carried by $\tau$, while keeping $c$ tight with respect to $\mathcal E$. Since $\tau$ is tight with respect to $\mathcal E$, it follows from conditions (1) and (3) of Definition \ref{tighttts} that every singular fiber of $\tau$ is a subarc of some $E_i \in \mathcal E$. Hence it suffices to show that any of the $c_i$ not contained in $\tau$ is contained in a possibly larger segment $c'_i$ of $c$ with endpoints $\partial c'_i$ on $\mathcal E$, and that $c'_i$ can be moved into $\tau$ by an isotopy relative to its endpoints, while keeping it tight with respect to $\mathcal E$. This will be done now by considering one-by-one each of the complementary components $\Delta$ of $\tau$ and moving $c$ off $\Delta$.

Since $\beta$ covers $\tau$, it follows from Lemma 2.4 of \cite{LM1} that for every complementary component $\Delta$ of $\tau$,  with sides ,say, $\delta_1, \delta_2, \delta_3$, there is for any $i = 1, 2, 3$ a subpath $\delta_i'$ of $\beta$ that runs parallel to all of  $\delta_i$. Let $I_1, I_2$ and $I_3$ be the singular $I$-fibers which contain the cusp points of $\Delta$. Let $\hat \Delta$ be the ``hexagon'' with sides alternatingly situated on one of the $\delta'_i$ or on one of the  $I_k$, and which contains 
$\Delta$.

Since $c$ is simple, any connected component $c'$ of $c \cap \hat \Delta$ must have endpoints on some $I_j$ and $I_k$. From the tightness of $c$ it follows $I_j \neq I_k$, since both of them are contained in curves from $\mathcal E$, by conditions (1) and (3) of Definition \ref{tighttts}.

Thus, as $\hat{\Delta}$ is simply connected,  $c'$ can be isotoped to an arc $c''$ that is contained in one of three connected components $\Delta^i$ of $\hat \Delta \smallsetminus \Delta$, and this can be done simultaneously with all such arcs $c'$ while keeping the curve $c$ simple, and also keeping it tight with respect to $\mathcal E$, by condition (2) of Definition \ref{tighttts}. But each $\Delta^i$ belongs to $\tau$, and we can assume that $c''$ is transverse to the $I$-fibers. 

It follows, after  performing all of these isotopies for any complementary component $\Delta$ of $\tau$, that $c$ is carried by $\tau$.

\end{proof}

\begin{remark} \label{correction} \rm
The statement of Proposition \ref{carriedcurves} includes that of Lemma 3.9 of \cite{LM1}. Unfortunately it seems that the proof given there is not quite correct; hence the above proof serves also as correction of the latter.

\end{remark}

\vskip8pt

\subsection{A distance criterion in the curve complex}
\hfill

\vskip10pt

Sequences of nested train tracks, as given in the previous  definition, occur already in \cite{MM}, 
Section 3.1, where they are used to derive  lower bounds for the  distance in the curve complex.  
Indeed, the  following  statement is a variant of their ``Basic  observation''.  A detailed proof is given in \cite{LM1}, stated there as Proposition 2.12 and Remark 2.13.

\vskip8pt

\begin{corollary}  \label{corlargedistance} 
For $n \geq 1$  let $\tau_{0} \supset  \tau_{1} \supset \ldots  \supset\tau_{n}$ be an $n$-tower  of derived train tracks in  $\Sigma$.  Assume that $\tau_0$ is a complete fat  train track.  Let $D$ be an   essential simple closed curve  carried by $\tau_n$, and let $E$ be an  essential simple closed   curve which has 
a wave with respect to  $\tau_0$.  Then one  has:
$$d_\mathcal C(D, E)\geq n+1$$

\end{corollary}

\vskip8pt

\subsection{Heegaard splittings}
\label{subsection-Heegaard-splittings}

\hfill

\vskip10pt

 Let $H$ be a handlebody of genus $g \geq 2$, and let $\Sigma =  \partial H$ denote its boundary surface.  The set  $\mathcal{D}(H)$  of isotopy classes of essential simple closed curves on $\Sigma$ that  bound a disk in $H$ is a subset of  $\mathcal{C}^0(\Sigma)$. It is  the vertex set  of what is called the {\it disk complex} of the  handlebody $H$,  contained as a subcomplex in $\mathcal{C}(\Sigma)$.  

\vskip5pt

Similarly, we consider complete decomposing systems, up to isotopy in  $\Sigma$, which bound disk systems in $H$, and denote the set of such  isotopy classes by $\mathcal{CDS}(H)$.

\vskip5pt

Any closed orientable $3$-manifold $M$ has a {\it Heegaard  splitting}, which is a decomposition 
of  $M$ along a surface $\Sigma$ into two genus $g$  handlebodies   $V$ and $W$,  so that 
$M = V \cup_{\Sigma} W$.   The genus of  the {\it Heegaard surface} $\Sigma$ is called the {\it genus} of the  Heegaard splitting. 

\vskip10pt

The {\it  distance} of a  Heegaard splitting  $M = V \cup_{\Sigma}  W$  is defined by
$$d(V, W)   = \min \{d_{\mathcal C}(D, E) \mid D\in \mathcal  {D}(V),  E \in \mathcal {D}(W) \}\, ,$$
where $d_{\mathcal C}$ denotes, as before, the distance in the  curve   complex $\mathcal{C}(\Sigma)$ 
(see \cite{He}).

\vskip10pt

\begin{remark} \label{otherdisks} \rm
Given a complete decomposing system 
$${\mathcal D} = \{D_1, ..., D_{3g-3}\} \in \mathcal{CDS}(V)$$ 
for a  handlebody $V$ of genus $g \geq 2$, then any other essential disk-bounding curve   
$D\in  \mathcal{D}(V)$ is either parallel to one of  $D_i $, or $D$   has a wave with respect to 
$\mathcal D\,$.

\end{remark}

\vskip10pt

A  complete decomposing system  ${\mathcal D} = \{D_1, ...,  D_{3g-3}\}  \subset \Sigma$ is 
said to   {\it have a wave with respect to a  second  complete  decomposing system 
${\mathcal E} \subset \Sigma$}  if some of the $D_{i}$ has a wave  with respect to  $\mathcal E$.
\vskip10pt

\begin{lemma}[\cite{He}, Lemma 1.3] \label{quotehempel}  
For every Heegaard splitting of a 3-manifold $M = V \cup_{\Sigma}  W$   there always  exists a 
pair of complete decomposing systems  ${\mathcal D} \in   \mathcal{CDS}(V)$ and  
${\mathcal E \in \mathcal{CDS}(W)}$ which have no waves with  respect   to each other.

\end{lemma}

\vskip10pt

The following is the main result of \cite{LM1}:

\vskip10pt

\begin{theorem}[Theorem 4.7 of \cite{LM1}] \label{HSdistance}
Let $M$ be an oriented $3$-manifold with a Heegaard splitting $M = V   \cup_{\Sigma} W$.   Consider  complete decomposing systems $\mathcal D \in  \mathcal{CDS}(V)$ and  
$\mathcal E \in \mathcal{CDS}(W)$ which do not have waves with  respect to each other.  
Let $\tau \subset \Sigma$ be a complete fat train track with    exceptional fibers 
$\mathcal{E}_{\tau} = \mathcal E$, and assume that $\mathcal D$ is   carried by $\tau_{n}$,  for some 
$n$-tower of derived train tracks   $\tau = \tau_{0} \supset  \tau_{1}  \supset \ldots \supset \tau_{n}$ with $n \geq 2$.  Then the distance of the  given Heegaard splitting  satisfies: $$d(V, W) \geq n$$

\end{theorem}

\vskip10pt

\section{Gregarious laminations}\label{gregarious}

\vskip10pt

A lamination $\mathcal L \subset \Sigma$ is called {\it minimal} if every leaf is dense in 
$\mathcal L$. A lamination $\mathcal L  \subset \Sigma$ is called  {\it filling} if   the components of 
$\Sigma \smallsetminus \mathcal L$ are simply connected. 
A {\it measured geodesic  lamination}  $({\mathcal L}, \mu) $ on $\Sigma$ is a lamination ${\mathcal L}$   together with a transverse measure  $\mu$  supported on $\mathcal L$ 
(see \cite{PH}). Such a measured lamination   $({\mathcal L}, \mu) $ is called {\it uniquely ergodic} if any   transverse measure supported on   $\mathcal L$ is  a multiple of $\mu$. As is common use, we denote the space of projective measured laminations $[{\mathcal L}, \mu]$ on $\Sigma$ by $\mathcal{PML}(\Sigma)$ (see~\cite{FLP} and   ~\cite{Ke}). It  comes with a natural measure class given by Thurston's  p.l.-structure of the  $(6g-7)$-dimensional sphere $\mathcal{PML}(\Sigma)$.

\vskip5pt

Recall that a train track $\tau$ is called maximal if every connected component of 
$\Sigma  \smallsetminus  \tau$ is a triangle.   

\vskip8pt

\begin{lemma} \label{minfilful}
The subset of $\mathcal{PML}(\Sigma)$  given by all minimal   laminations has full measure in  
$\mathcal{PML}(\Sigma)$.

\end{lemma}

\vskip4pt

\begin{proof} 
It is well known (see ~\cite{Ma}, ~\cite{Ve}) that the set of uniquely ergodic laminations has full  measure in $\mathcal{PML}(\Sigma)$.  We only need to consider geodesic laminations $\mathcal L$ that are  given as the support of some transverse measure $\mu$ carried by $\mathcal L$.  Since every such  lamination which is not minimal  is also  non-uniquely ergodic, the  set of minimal filling laminations contains  the uniquely ergodic   ones, which shows the desired conclusion.

\end{proof}

\vskip8pt

Any uniquely ergodic measured lamination $({\mathcal L}, \mu)$ has the property that the 
geodesic lamination ${\mathcal L}$ determines   the corresponding  projective class 
$[{\mathcal L}, \mu] \in \mathcal{PML}(\Sigma)$.  This justifies a    certain amount of sloppyness in suppressing the difference between  laminations and projective classes of measured laminations.

\vskip10pt

\begin{definition} \rm 
Let $\tau \subset \Sigma$ be a train track. We will use the   following  notation:

\begin{enumerate}
\item[(1)] $\mathcal{PML}(\tau)$ is the set of projective measured  laminations carried by $\tau$. 

\vskip5pt

\item[(2)] $\bf{P}(\tau) \subset \mathcal{PML}(\tau)$ is the subset  of all  projective measured laminations which are carried by $\tau$  and have   positive weights on every edge of $\tau$.  
Such a  lamination will  be   called $\tau$-positive. 

\vskip5pt

\item[(3)]  $\bf {M}(\tau) \subset \mathcal{PML}(\tau)$ denotes the  subset  given by all minimal laminations which are carried by $\tau$.  Set   $\bf  {MP}(\tau) = \bf M(\tau) \cap \bf P(\tau)$.

\vskip5pt

\item[(4)]   $\bf {A}(\tau) \subset \mathcal{PML}(\tau)$ denotes the  subset  given by all {\it arational} laminations  (i.e. not given by  a collection of  simple closed curves) which are   carried by $\tau$. 
Set   $\bf {AMP}(\tau) =  \bf A(\tau) \cap \bf M(\tau) \cap \bf P(\tau)$.

\end{enumerate}

\end{definition}

\vskip8pt

To be specific, let us observe here that an element $[\mathcal{L}, \mu]$ of ${\bf M}(\tau)$ does 
not belong to $\bf {A}(\tau)$  if and  only if the support of $\mu$  consists of a  single closed curve. 
If  $[\mathcal{L},  \mu] \in{\bf M}(\tau) \cap \bf {A}(\tau)$,  then the support  of $\mu$ is {\it totally  arational}:  it does not  contain any closed leaf.

\vskip10pt

\begin{lemma}\label{MPfull}
Let  $\tau \subset \Sigma$ be a  maximal train track. Then one has:

\begin{enumerate}

\item[(1)] The set $\bf{P}(\tau)$  is open in $\mathcal{PML}(\Sigma)$.

\vskip4pt

\item[(2)] The set $\mathcal{PML}(\tau)$  has positive measure in   $\mathcal{PML}(\Sigma)$.

\vskip4pt

\item[(3)] The set $\bf {AMP}(\tau)$ has full measure in   $\mathcal{PML}(\tau)$.

\end{enumerate}

\end{lemma}

\vskip10pt

\begin{proof} 
Since $\tau$  is maximal, it has only triangles as complementary   components. Hence the set 
$\mathcal{PML}(\tau)$ is a top  dimensional  cell in $\mathcal{PML}(\Sigma)$ and hence  is has 
positive  measure.  Furthermore, Lemma \ref{minfilful} implies that   $\bf{M}(\tau)$ has  full measure  in 
$\mathcal{PML}(\tau)$.  On the other hand,  $\bf{P}(\tau)$ is   precisely the interior of this top 
dimensional  cell, and hence  it  is open and has full measure in   $\mathcal{PML}(\tau)$.
Thus $\bf {MP}(\tau) = \bf{M}(\tau)\cap \bf{P}(\tau)$ has full   measure in $\mathcal{PML}(\tau)$. 
But the set of rational laminations is  countable and hence of   measure $0$.  Thus $\bf {AMP}(\tau)$  
has also  full measure in $\mathcal{PML}(\tau)$.

\end{proof}

\vskip6pt

\begin{definition}\label{gregarious} \rm  Let $\tau$ be a train track on $\Sigma$.
\begin{enumerate}
\item[(a)]  A lamination $\mathcal L$ carried by  $\tau$ is called  {\it  gregarious}\footnote{\, The authors would  like to thank Wendy Sandler  for suggesting this inspiring terminology.} with  respect to $\tau$ 
if  the train track $\tau$ can be derived with respect to $\mathcal L$ (compare Definition \ref{derivedtt}).

\vskip4pt

\item[(b)] The subset of $\mathcal{PML}(\tau)$ defined by all   gregarious laminations with respect to 
$\tau$ will be denoted by $\bf{G}(\tau)$.  Set  $\bf {GP}(\tau) =  \bf  G(\tau) \cap \bf P(\tau)$.

\end{enumerate}

\end{definition}

\vskip8pt

\begin{lemma} \label{minimalgregarious} 
Let $\tau$ be a maximal train  track  on $\Sigma$, and let $\mathcal {L}$   be an arational minimal  
$\tau$-positive lamination. Then $\mathcal  L$  is gregarious  with respect to $\tau$: 
$$\bf{AMP}(\tau) \subset \bf{GP}(\tau)$$

\end{lemma}

\begin{proof}
This is a direct consequence of  Definition \ref{gregarious}: Since   $\mathcal {L}$ is $\tau$-positive
it covers $\tau$, and since it is arational, no two zippers can ever   meet, for arbitrary long  unzipping  paths. Since $\mathcal {L}$ is minimal, every  unzipping path will  eventually intersect every transverse 
$I$-fiber  which is met by $\mathcal {L}$.   But $\mathcal {L}$ is minimal and covers  $\tau$, so that the unzipping paths will eventually become  complete.

\end{proof}

The next lemma states that the inclusion from the previous lemma, though an equality ``in  measure'' 
(by Lemma \ref{MPfull} (3)),  is proper,  since ${\bf AMP}(\tau)$ is disjoint from the dense set of laminations  supported by a single closed curve.

\vskip8pt

\begin{lemma} \label{gregariousopen}  
Let $\tau \subset \Sigma$ be a maximal train track. Then the  set   $\bf {GP}(\tau)$ is  open  in 
$\mathcal{PML}(\tau)$.

\end{lemma}

\vskip5pt

\begin{proof} 
Let $\mathcal L$ be a $\tau$-positive lamination which   is gregarious with respect  to $\tau$ and let 
$\tau' = \tau(\mathcal{L})$ be the  derived train track.   Since $\tau$ is maximal, it follows that $\tau'$ 
also has only triangles as complementary  components. An elementary Euler characteristic count 
shows that  $\mathcal{PML}(\tau')$ is a  subcell of maximal dimension in  $\mathcal{PML}(\tau)$, 
and  thus any open subset of $\mathcal{PML}(\tau')$  is open in $\mathcal{PML}(\tau)$.

Every element $[{\mathcal L}, \mu]$ of 
$\bf{GP}(\tau) \subset  \mathcal{PML}(\tau') \subset  \mathcal{PML}(\tau) \subset \mathcal{PML}(\Sigma)$ defines a set of   weights on the edges of $\tau'$, where we impose the additional
condition that   the sum of these weights is one: Otherwise the weights would only be   determined up to a scalar factor.

Perturbing the weights on the edges of $\tau'$  slightly determines   an open neighborhood of 
$[{\mathcal L}, \mu] \in \bf{GP}(\tau)$ in   $\mathcal{PML}(\tau')$.  This  neighborhood is also open 
in $\mathcal{PML}(\tau)$, and  (by Lemma \ref{MPfull} (1)\,) even  in  $\mathcal{PML}(\Sigma)$. Since any lamination in this neighborhood is  carried by  $\tau'$,  the neighorhood consists  entirely of laminations in $\bf{GP}(\tau)$. Thus $\bf{GP}(\tau)$ is   open in  $\mathcal{PML}(\tau)$. 

\end{proof}

\vskip5pt

\begin{lemma} \label{zippingpaths}
Let $\tau \subset \Sigma$ be a maximal train track.  Then the  set  $\bf  {GP}(\tau)$ has full measure 
in  $\mathcal{PML}(\tau)$.
\end{lemma}

\vskip5pt

\begin{proof}
This is a direct consequence of Lemma \ref{minimalgregarious} and   
Lemma \ref{MPfull} (3).

\end{proof}

\vskip5pt

We will now consider the the subset $\mathcal{PML}_{0}(\Sigma)$ of   $\mathcal{PML}(\Sigma)$ 
which consists  of projective measured laminations $[\mathcal{L}, \mu]$ such that $\mathcal L$ is
a single closed curve. Since a simple closed curve carries  up to scalar multiples only one  transverse measure, there is a canonical identification  
$\mathcal{PML}_{0}(\Sigma) = \mathcal{C}^{0}(\Sigma)$.

In  analogy to $\mathcal{PML}_{0}(\Sigma)$, we introduce   
$\mathcal{PML}_{0}(\tau) = \mathcal{PML}(\tau) \cap    \mathcal{PML}_{0}(\Sigma)$,  
$\bf{GP}_{0}(\tau) = \bf{GP}(\tau) \cap \mathcal{PML}_{0}(\Sigma)$,   
etc. Of course, all these newly introduced sets are countable, and they    are  dense in  their ``parent'' 
set, if the latter is open in $\mathcal{PML}(\Sigma)$.

\vskip8pt

\begin{lemma}\label{countablymanytts} 
Let $\tau \subset \Sigma$ be a maximal train track.  Then there  is  a  countable  family of  maximal
train tracks  $\tau_{1}, \tau_{2}, \ldots$, each of which is  derived from $\tau$,   such  that  
${\bf GP}(\tau)$ is equal to the union of all ${\bf{P}}(\tau_{i})$.

\end{lemma}

\begin{proof}
Since $\mathcal{PML}_{0}(\Sigma)$ is dense in $\mathcal{PML}(\Sigma)$  
and ${\bf GP}(\tau) \subset {\bf P}(\tau)$ is open in  $\mathcal{PML}(\Sigma)$ (by Lemma 
\ref{gregariousopen}  and  Lemma \ref{MPfull} (1)), it follows that the countable set $\bf{GP}_{0}(\tau)$ is dense in  $\bf{GP}(\tau)$.  Notice that  for  each  element $D_{i}$ of   $\bf{GP}_{0}(\tau)$  there is a 
maximal train  track $\tau_{i}$  derived from $\tau$ which  carries  $D_{i}$, and  that  ${\bf{P}}(\tau_{i})$
is an open  neighborhood of  $D_{i}$ in $\mathcal{PML}(\tau)$, again  by Lemma \ref{gregariousopen}.
Hence the union of all ${\bf{P}}(\tau_{i})$  contains all of ${\bf  GP}(\tau)$, and hence is equal to the latter.

\end{proof}

The notion of  gregariousness can be  strengthened  further:  

\vskip8pt

\begin{definition} \label{ngregarious}\rm
Let $\tau \subset \Sigma$ be a train track, and let $\mathcal L$ be  a  lamination carried by $\tau$. 
We say that $\mathcal L$ is {\it  $n$-gregarious} with respect to    $\tau$ if $\tau$ can be derived 
$n$ times with respect to  $\mathcal L$,  i.e. there exists a tower  
$$\tau = \tau_{0} \supset \tau_{1} \supset \ldots \supset \tau_{n}$$ 
of derived train tracks  with respect to $\mathcal L$. In particular  $\mathcal L$ is carried by $\tau_{n}$. 
We denote the subset  of  $\mathcal{PML}(\tau)$  given by all   $n$-gregarious laminations by 
${\bf G}^n(\tau)$, and  define   $\bf{G}^{n}\bf{P}(\tau) = {\bf   G}^n(\tau) \cap {\bf P}(\tau)$.

\end{definition}

\vskip5pt

\begin{proposition} \label{doublezippingpaths}
Let $\tau \subset \Sigma$ be a maximal train track. Then for  all   $n\geq 1$ the set   
$\bf{G}^{n}\bf{P}(\tau)$ of $n$-gregarious laminations which are   $\tau$-positive is open  in  
$\mathcal{PML}(\tau)$  and in $\mathcal{PML}(\Sigma)$.
\end{proposition}

\vskip5pt

\begin{proof}
By Lemma \ref{MPfull} (1) it suffices to prove openess in  $\mathcal{PML}(\tau)$. Assume by induction that  $\bf{G}^{n-1}\bf{P}(\tau)$ is open  in   $\mathcal{PML}(\tau)$, and that there is a countable family of maximal  train tracks $\tau_{i}$ that are  obtained from $\tau$ by deriving $n-1$ times, such that 
$\bf{G}^{n-1}\bf{P}(\tau)$ is equal to the union of all ${\bf P}(\tau_{i})$. Thus $\bf{G}^{n}\bf{P}(\tau)$ is equal to the union of all ${\bf GP}(\tau_{i})$. 

Now apply Lemma \ref{countablymanytts} to  each of the $\tau_{i}$ to get a countable family of 
maximal train tracks  $\tau_{i, j}$ derived from $\tau_{i}$, such that  ${\bf GP}(\tau_{i})$ is equal to the union of all ${\bf P}(\tau_{i,j})$, for fixed $i$.  It follows that  $\bf{G}^{n}\bf{P}(\tau)$ is equal to the union of the ${\bf P}(\tau_{i, j})$, for all $i$ and $j$. 

From Lemma \ref{gregariousopen} we obtain that every ${\bf P}(\tau_{i, j})$ is open in 
$\mathcal{PML}(\tau_{i,j})$.  Since  all $\tau_{i,j}$ are maximal, the set $\mathcal{PML}(\tau_{i,j})$ 
is a cell of  maximal dimension in $\mathcal{PML}(\tau)$. It follows directly that  every ${\bf P}(\tau_{i, j})$
is open in $\mathcal{PML}(\tau)$. Hence their union $\bf{G}^{n}\bf{P}(\tau)$ is also open in 
$\mathcal{PML}(\tau)$. This completes the  induction and hence the proof.

\end{proof}

\vskip0pt

\begin{proposition} \label{nfullmeasure}
Let $\tau \subset \Sigma$ be a maximal train track. Then for all $n\geq 1$ the set 
$\bf{G}^{n}\bf{P}(\tau)$ of  $n$-gregarious laminations which are $\tau$-positive 
is of full measure in $\mathcal{PML}(\tau)$.

\end{proposition}

\vskip0pt

\begin{proof}
By Lemma \ref{MPfull} (3) we know that $\bf{G^{0}}\bf{P}(\tau) = {P}(\tau)$ has full measure  in 
$\mathcal{PML}(\tau)$.  Hence the claim follows by induction if one proves that $\bf{G}^{k+1}\bf{P}(\tau)$ has full  measure  in    $\bf{G}^{k}\bf{P}(\tau)$.

Let us first recall that $\bf{G}^{k}\bf{P}(\tau)$ is the countable union of sets ${\bf P}(\tau_{i})$, where each $\tau_{i}$ is a maximal train track obtained  from  $\tau$ by deriving  $k$ times. This has been shown by induction in the proof of Proposition  \ref{doublezippingpaths}.

For any $\tau_{i}$ the set ${\bf{GP}}(\tau_{i})$ has  full  measure in  
${\bf P}(\tau_{i}) \subset \mathcal{PML}(\tau_{i})$ (by Lemma \ref{zippingpaths}). Hence the  union 
of all ${\bf{GP}}(\tau_{i})$ has full measure in the union of all ${\bf P}(\tau_{i})$.   But  the  union of all 
${\bf{GP}}(\tau_{i})$  is (by   definition of the $\tau_{i}$) equal to $\bf{G}^{k+1}\bf{P}(\tau)$, 
while the union of all ${\bf P}(\tau_{i})$ is precisely   $\bf{G}^{k}\bf{P}(\tau)$. 
This proves the inductive  step.
\end{proof}

\begin{proposition}
\label{ngregariousrem}
For any complete fat train track $\tau$ in $\Sigma$ the set ${\bf G^n_{0}}(\tau)$ has  closure 
$\overline{{\bf G^n_{0}}(\tau)}$ in $\mathcal{PML}(\tau)$ which is of full measure. The complement of 
${\bf G^n_{0}}(\tau)$  in $\mathcal{PML}_{0}(\tau)$ has closure 
$\overline{\mathcal{PML}_{0}(\tau) \smallsetminus {\bf G^n_{0}}(\tau)}$ in $\mathcal{PML}(\tau)$
which is of measure 0.
\end{proposition}

\begin{proof}
The set $\mathcal{PML}_{0}(\Sigma)$ is dense in  $\mathcal{PML}(\Sigma)$, and thus,  since  
$\bf{G}^{n}\bf{P}(\tau)$  is  open in  $\mathcal{PML}(\Sigma)$ (by Proposition \ref{doublezippingpaths}),
it follows that   $\bf{G}^{n}\bf{P}_{0}(\tau) =  \mathcal{PML}_{0}(\Sigma) \cap  \bf{G}^{n}\bf{P}(\tau)$ 
is dense  in $\bf{G}^{n}\bf{P}(\tau)$. But  by Proposition  \ref{nfullmeasure}  the set $\bf{G}^{n}\bf{P}(\tau)$ 
is of full measure in $\mathcal{PML}(\tau)$. Since $\bf{G}^{n}\bf{P}_{0}(\tau)$ is a subset of 
$\bf{G}^{n}_{0}(\tau)$,  it follows that $\overline{{\bf G^n_{0}}(\tau)}$ has full measure
in $\mathcal{PML}(\tau)$.

The complementary set $\mathcal{PML}_{0}(\tau) \smallsetminus {\bf G^n_{0}}(\tau)$ is contained in 
$\mathcal{PML}_{0}(\tau) \smallsetminus {\bf G^{n}P_{0}}(\tau)$,  and hence in 
$\mathcal{PML}(\tau) \smallsetminus {\bf G^{n}P}(\tau)$, since
${\bf G^{n}P}(\tau) \cap \mathcal{PML}_{0}(\tau) = {\bf G^{n}P_{0}}(\tau)$. By Propositions 
 \ref{doublezippingpaths} and \ref{nfullmeasure}, the set 
$\mathcal{PML}(\tau) \smallsetminus {\bf G^{n}P}(\tau)$ is a closed set of measure 0 in  
$\mathcal{PML}(\tau)$.

\end{proof}

\vskip0pt

\section{Genericity}\label{generic}

\vskip10pt

Classically, a  subset of a countable set is called ``generic'' if  its  complement  is finite.  This notion, however, often doesn't  capture  the geometry of the given set-up.

\vskip4pt

For example, consider the countable set $S$ of points in the unit square $I^{2}$ which have rational 
coordinates. The subset of $S$ which lies in the interior of $I^{2}$   has infinite complement, but everyone will agree that a  ``generic'' point  of $S$ will lie in the interior of  $I^{2}$ and  not on its boundary.

\vskip4pt

In order to address the above problem we propose the following  more subtle definition for  genericity:

\vskip6pt

\begin{definition}\label{generic} \rm
Let $X$ be a topological space,  provided with a Borel measure  $\mu$. Let $Y \subset X$   be a (possibly countable) subset, which is a disjoint  union   $Y =  A \overset {\cdot}{\cup}  B$.   The set $A$ is called {\it generic in $Y$} (or simply {\it  generic},  if  $Y = X$) if the closure  $\bar A$ of $A$ has measure  $\mu(\bar A) > 0$, and  the closure    $\bar B$ of $B$ has   measure $\mu(\bar B) = 0$.
To be specific, both closures $\bar A$ and $\bar B$ are taken in $X$.

\end{definition}

\vskip4pt

\begin{remark}\rm Notice that it follows from Definition \ref{generic} that in the above setting, whether or not a set $A$ is generic, depends  only on the measure class of $\mu$, and not on the measure $\mu$ 
itself.

\end{remark}

Notice that in this definition the sets $\bar  A$  and $\bar B$ may well not be disjoint, although $A$ 
and $B$ are  assumed to be disjoint.  Note also, that this definition of genericity  extends to sets $Y$ that are not embedded but are just mapped to $X$, by a properly chosen ``natural'' map.   It is  important 
to remember that  every statement about genericity always depends on  a previous choice  of a
measure. This  choice is,  formally speaking, arbitrary, and thus can at best be  natural.

\vskip8pt

The following is an immediate  consequence of the above definition:

\vskip8pt

\begin{lemma} \label{equivalentgeneric}
Given sets  $X, Y, A$ and $\bar A$ as in Definition \ref{generic}, and let $\bar Y$ be the closure of $Y$ in $X$. Then $A$  is  generic in $Y$ if and only if one of the following two  equivalent  conditions is satisfied:

\begin{enumerate}
\item 
The set  $\bar A$ contains a set $Z$ which is open in  $\bar Y$ and  of full measure  
$\mu(Z) = \mu(\bar Y) > 0$, and which is disjoint from $Y  \smallsetminus A$.

\item The set $\bar Y$ has measure  $\mu(\bar Y) > 0$, and $Y \smallsetminus A$  is contained 
in a closed set of measure 0.

\end{enumerate}

\qed

\end{lemma}

\vskip5pt

\begin{remark} \label{reformulation} \rm
As a direct consequence of Definition \ref{generic} and its reformulations in Lemma  
\ref{equivalentgeneric} we  obtain the following:

\begin{enumerate}

\item[(a)] Arbitrary unions and finite  intersections of sets $A_{i}$  that are generic in a common set 
$Y  \subset X$ are  again generic in $Y$.  

\item[(b)] For any sets $A \subset A' \subset Y \subset X$, if  $A$ is generic in $Y$,  then so is $A'$.

\end{enumerate}

\end{remark}

\vskip5pt

The situation becomes more  complicated if one also  varies the  set  $Y$.  

\vskip5pt

\begin{proposition} \label{unionofgenerics}
Let $X, Y$ and $A$ be as in Definition \ref{generic}. Assume that $Y$  contains a countable union of 
subsets $Y_{i}$, and define $A_{i} = Y_{i} \cap A$.  Denote, as before, by $\bar Y$ and $\bar Y_i$ the closures  (in $X$) of $Y$ and $Y_i$ respectively, and assume furthermore that:

\begin{enumerate}

 \item[(a)]   $Y_{i} = \bar Y_{i}  \cap Y$, 
 
\item[(b)] $\mu(\bar Y \smallsetminus \cup \bar Y_{i}) = 0$, and

\item[(c)] every $\bar Y_{i}$ contains some set   $\overset{\circ}{Y_{i}}$ that is open in $\bar Y$ 
and has full  measure in $\bar Y_{i}$.

\end{enumerate}

\noindent If every $A_{i}$ is generic in $Y_{i}$, then $A$  is generic in $Y$.

\end{proposition}

\vskip5pt

\begin{proof}
By assumption every $A_{i}$ is generic  in $Y_{i}$.  Hence Lemma  \ref{equivalentgeneric} (1) gives sets
$Z'_{i} \subset \bar A_{i}$ that are open in $\bar   Y_{i}$, are of full positive measure  in  $\bar Y_{i}$, and satisfy  $Z'_{i} \cap Y_{i} \subset A_{i}$. Define  $Z_{i} = Z'_{i} \cap \overset{\circ}{Y_{i}}$, and observe that  the  $Z_{i}$ are still of full  positive measure in  $\bar Y_{i}$, and in addition  they are open in $\bar Y$. Their union  $Z = \cup Z_{i}$ is  an open set (in $\bar Y$) of  positive  measure contained in $\bar A$, so that  $\mu(\bar A) > 0$.

The complementary set $\bar Y \smallsetminus Z$ is closed in $\bar Y$ and hence in $X$, 
and  it  contains $Y \smallsetminus A$,  since for any index $i$ one has $Z_{i} \subset \bar Y_{i}$, and 
$Z_{i} \cap Y \subset Z'_{i} \cap Y = Z'_{i} \cap \bar Y_{i} \cap Y =  Z'_{i} \cap Y_{i} \subset A_{i}$.
But $\bar Y \smallsetminus Z$  is contained in 

\vskip5pt

 \centerline{$(\bar Y  \smallsetminus \cup \bar Y_{i}) \, \, \cup \, \, (\,\cup \,  (\bar Y_{i} 
 \smallsetminus  Z_{i})  ) \, ,$} \hfill

\noindent 
which is a  countable union of measure 0 sets. Thus $\bar Y \smallsetminus Z$ is of measure 0, which implies that   $\bar{Y \smallsetminus A}$ is of measure 0. This proves that $A$ is   generic in $Y$.

\end{proof}

\vskip7pt

The following proposition will not be used below, but we believe it  can be a useful tool in other contexts.

\begin{proposition} \label{genericitycriterion}
Let $X$  and $Y$ be as in Definition \ref{generic}.  Assume that  the set $A \subset Y$ is given as (a not necessarily disjoint) union  of countably  many subsets $A_{i}$, and for each $A_{i}$ there are 
sets $X_{i}$ and $Z_{i}$ such that for any index $i$ the following  holds:

\begin{enumerate}
\item[(1)] $Z_{i} \subset \bar A_{i} \subset X_{i} \subset X$
\vskip4pt
\item[(2)] $X_{i}$ is closed, and $Z_{i}$ is open in  $X$.
 \vskip4pt   
\item[(3)] The union of all $X_{i}$  contains $Y$.
 \vskip4pt
\item[(4)] $Z_{i} \cap Y \subset A_{i}$.
\vskip4pt
\item[(5)] $Z_{i}$ has full measure in $X_{i}$.
\vskip4pt
\item[(6)] Some $X_{j}$ has positive measure in $X$.
\vskip4pt
\item[(7)] $\bar{\cup X_{i}} \smallsetminus \cup X_{i}$ has measure 0.

\end{enumerate}

\vskip4pt

Then $A$ is generic in $Y$.

\end{proposition}

\vskip10pt

\begin{proof}
By (1) $Z_{i}$ is contained in $\bar A$, and by (5) and (6) some  $Z_{i}$  has positive measure.  Thus 
$\bar A$ has positive measure. On the other hand, one has:
$$\bar{\underset{i}{\cup}  X_{i}} \smallsetminus \underset{i}{\cup} 
Z_{i} \, \,\subset \, \,  \bar{\underset{i}{\cup}  X_{i}} 
\smallsetminus \underset{i}{\cup}  X_{i}  \, \, \cup \, \,  
\underset{i}{\cup}  X_{i} \smallsetminus \underset{i}{\cup} Z_{i} $$
and
$$\underset{i}{\cup}  X_{i} \smallsetminus \underset{i}{\cup} Z_{i}  
\, \,\subset \, \,  \underset{i}{\cup} (X_{i} \smallsetminus Z_{i})$$
The set  $\underset{i}{\cup} (X_{i} \smallsetminus Z_{i})$ has  measure 0, by (5) and by the countability of the index set.   Similarly, the set 
$\bar{\underset{i}{\cup}  X_{i}} \smallsetminus \underset{i}{\cup}   X_{i}$  has measure 0, by assumtion (7). 
The set $\bar{\underset{i}{\cup}  X_{i}} \smallsetminus \underset{i}{\cup} Z_{i}$ is closed, by (2).
But by (3) and (4) the set $Y \smallsetminus A$ is contained in  
$\bar{\underset{i}{\cup}  X_{i}} \smallsetminus \underset{i}{\cup}  Z_{i}$,  so that the closure  
$\bar{Y \smallsetminus A}$ must have measure 0.  Thus $A$ is generic in $Y$.

\end{proof}

\vskip5pt

\section{Genericity of large distance in the curve complex}
\label{curves}

\vskip10pt

For any handlebody $H$ with boundary $\partial H = \Sigma$ and any integer $n \in \naturals$ 
we say that a curve $c \in \mathcal{C}^{0}(\Sigma) = \mathcal{PML}(\Sigma)$ is {\it $n$-gregarious 
with respect to $H$}, if $\,c\,$ is  $n$-gregarious with respect to  some complete fat train track $\tau$ 
with exceptional fibers $\mathcal E_{\tau}$ in $\mathcal{CDS}(H)$  (compare with Definition 
\ref{ngregarious}). In the terminology of Definition \ref{gregariousnew}, this is equivalent to stating that 
$c$ is $n$-gregarious with respect to any complete decomposing system 
$\mathcal{E} \in \mathcal{CDS}(H)$.

\vskip8pt

\begin{proposition}
\label{morelongdistance}
For any integer $n \geq 1$ and any handlebody $H$ with boundary surface  $\partial H = \Sigma$ the set ${\bf G^n_0}(H)$ of $n$-gregarious curves  $c$ with respect to $H$ is generic in $\mathcal{C}^{0}(\Sigma)$.
\end{proposition}

\vskip4pt

\begin{proof}
We will use Proposition \ref{unionofgenerics}, with 
$X = \mathcal{PML}(\Sigma), Y = \mathcal C^{0}(\Sigma) =  \mathcal{PML}_{0}(\Sigma), 
A =  {\bf G^n_0}(H)$,  and $Y_{i} =  \mathcal{PML}_{0}(\tau_{i})$,  where $\tau_{i}$ is any complete fat  train track with  $\mathcal  E_{\tau_{i}} \in \mathcal{CDS}(H)$.  Note that the set $\mathcal{CDS}(H)$ of 
isotopy classes of complete decomposing systems  in $H$  is countable, and that for each complete decomposing system   $\mathcal E$ there are only  countably many fat train tracks  $\tau$ with 
$\mathcal{E}_{\tau} = \mathcal E$  (up to isotopy of the pair $(\Sigma, \mathcal E)$).  Note also 
that  $\mathcal{PML}_{0}(\tau_{i})$ is  dense in  $\mathcal{PML}(\tau_{i})$, which is closed in 
$\mathcal{PML}(\Sigma)$,  so that $\bar Y_i = \mathcal{PML}(\tau_{i})$ and $Y_{i} =  \bar Y_{i} \cap Y$ holds. 

Furthermore, we know from Lemma \ref{MPfull} that the set 
${\bf P}(\tau_{i}) \subset \mathcal{PML}(\tau_{i})$  is open and full measure in  $\mathcal{PML}(\tau_{i})$.  Since $\tau_{i}$ is  maximal, the set  $\mathcal{PML}(\tau_{i})$ is a cell of  maximal dimension in 
$\mathcal{PML}(\Sigma)$, so that ${\bf P}(\tau_{i})$  is open in $\bar Y = \mathcal{PML}(\Sigma)$. 
Thus we can define   $\overset{\circ}{Y_{i}} = {\bf P}(\tau_{i})$.

We now consider the set $\bar Y \smallsetminus \cup \bar Y_{i}\,$: It   consists of all laminations 
$\mathcal L$  which are not carried by any complete fat  train track with  exceptional fibers in 
$\mathcal{CDS}(H)$. Thus, by Lemma \ref{fatcarries}, $\mathcal L$ must have a wave with respect 
to any complete  decomposing  system that bounds disks in $H$.   But the set of such laminations 
$\mathcal L$ is precisely the set  $\mathcal R  \subset \mathcal{PML}(\Sigma)$  for which Kerckhoff shows  $\mu(\mathcal R) = 0$, in his  proof that the limit set of the  handlebody group has measure 
$0$ (see  \cite{Ke}).

We can now apply Proposition \ref{ngregariousrem}  to each of the $\tau_{i}\, $: It states precisely that 
$A_{i} = {\bf G^{n}_{0}}(\tau_{i}) =  \mathcal{PML}_{0}(\Sigma) \cap  {\bf G}^{n}(\tau_{i}) = 
\mathcal{PML}_{0}(\tau_{i}) \cap  {\bf G}^{n}(H) $   is generic in 
$Y_{i} =  \mathcal{PML}_{0}(\tau_{i}) $. Thus  Proposition \ref{unionofgenerics} gives 
the desired conclusion.

\end{proof}

\vskip0pt

Denote by  $\mathcal C^{0}_{n}(H)  \subset \mathcal{PML}_{0}(\Sigma)$
the set of essential simple closed  curves  $D \in  \mathcal{PML}_{0}(\Sigma)$ which satisfy 
$d_{\mathcal C}(D, E) \geq n$  for any $E \in \mathcal D(H) \, $:
$$\mathcal C^{0}_{n}(H)  = \{D \in  \mathcal{PML}_{0}(\Sigma) \mid d_{\mathcal C}(D, 
\mathcal D(H) ) \geq n\}$$

\vskip20pt

Recall that any complete fat train track $\tau$ on the  surface  $\Sigma$ defines  a handlebody 
$H = H(\tau)$ with  boundary $\partial H = \Sigma$  by the  condition 
$\mathcal E_{\tau} \in \mathcal{CDS}(H)$, i.e. all $E_{i} \in   \mathcal E_{\tau}$ bound disks in $H$.  

\vskip5pt

\begin{theorem} \label{longdistancefromH} 
For every handlebody $H$ with $\partial H = \Sigma$  the set $\mathcal C^{0}_{n}(H)$ is generic in the set  $\mathcal  C^{0}(\Sigma)$.
\end{theorem}

\vskip5pt

\begin{proof}

For every complete fat train track $\tau$ on $\Sigma$ it follows from Remark \ref{otherdisks} that 
every disk in $\mathcal{D}(H(\tau))$ has a wave with  respect to $\tau$. Thus it follows from 
Corollary \ref{corlargedistance} that $\mathcal  C^{0}_{n}(H)$ contains  ${\bf G^n_0}(H)$.  The latter is generic in  $\mathcal{C}^{0}(\Sigma)$,  by Proposition  \ref{morelongdistance}. Thus an application of statement (b)  of Remark \ref{reformulation}  finishes the proof.

\end{proof}

\vskip5pt

\begin{corollary} \label{longdistancefromc}
For any essential simple closed curve $c$ on $\Sigma$, the set  $\mathcal C^{0}_{n}(c)$ of all 
essential simple closed  curves $k$ on $\Sigma$ with distance  
$$d_{\mathcal  C}(k, c) \geq n$$
is generic in  the set  $\mathcal C^{0}(\Sigma)$ of all essential simple closed  curves on  $\Sigma$.

\end{corollary}

\vskip8pt

\begin{proof}
Consider any handlebody $H$ which contains a disk with boundary  curve $c$, and observe that 
$c \in \mathcal{D}(H)$ implies 
$$\mathcal C^{0}_{n}(H) \subset \mathcal C^{0}_{n}(c)  \, .$$
Thus statement (b)  of Remark \ref{reformulation} gives directly the stated claim.

\end{proof}

\vskip8pt

\section{Intersection and Dehn twists}

\vskip8pt

Let $\Sigma$ be an orientable surface of genus  $g \geq 2$, and let $\mathcal {E}$ be a complete decomposing system  for $H$.  Let $k$ be an essential simple closed curve which is tight (see  
subsection \ref{fattraintracks}) with respect to the complete decomposing system $\mathcal E$ on 
$\Sigma$.  The number   of  intersection  points of  $k$  with $\mathcal E$ is called the  {\it $\mathcal E$-length}  of $k$  and is denoted by $\mid k \mid_{\mathcal E}$.   The same definition and notation will be  used for a simple arc $\alpha$ instead of $k$. However, in this case we always require that $\partial \alpha$ is contained in $\mathcal E$, and we count the two points of $\partial \alpha$ as  intersection points when we determine the $\mathcal E$-length of $\alpha$.

Two tight simple arcs   on $\Sigma$ are called  {\it parallel} (with respect to  $\mathcal E$) if they are isotopic to each other  via an isotopy of the  pair $(\Sigma, \mathcal{E})$.  In this  case it follows (but this is not equivalent !) that the arcs can be oriented so that their intersections  with 
$\mathcal E$ 
occur  at precisely the same sequence of curves $E_j \in \mathcal E$, and from the  same direction.

We also need to specify what we mean below by a {\it arc $c'$ on a closed curve $c$}:  Such an arc $c'$ is not  necessarily a subarc of $c$,
it can also be the image of a  path which is  immersed in $c$ but not embedded in $c$. In particular, $c'$ can wind around $c$ several times.
 
Let $c$ and $k$ be distinct essential simple closed curves on $\Sigma$ that  are tight with respect 
to $\mathcal E$, and let $P \in c\ \cap\ k$ be some intersection point.   We now consider maximal 
parallel arcs $\alpha$ on $k$ and  $\alpha'$ on $c$ such that $P$ is contained in 
$\hat\alpha$ and in $\hat\alpha'$, 
where
$\hat\alpha$ and $\hat\alpha'$ denote the extensions of $\alpha$ on $k$ and $\alpha'$ on $c$, across the pair-of-pants adjacent to the curves of ${ \mathcal E}$ that contain an endpoint of $\alpha$ and $\alpha'$. We call $\alpha$ the {\it intersection arc of $P$ on $k$} (and $\alpha'$ the {\it intersection arc of $P$ on $c$}). The length of either is called the {\it intersection length} of $k$ and $c$ at $P$ and denoted by $|P|_\mathcal E$, i.e.:
$$\mid P \mid_{\mathcal E} \, \, =\, \,  \mid \alpha \mid_{\mathcal E} 
\, \, \,=\, \,  \mid \alpha' \mid_{\mathcal E}\,$$
Note that, as every pair-of-pants complementary to $\mathcal E$ has precisely three boundary curves, the intersection arcs on $c$ and $k$ are well defined by $P$. 
Furthermore, one has always $|P|_\mathcal E \geq 1$ unless $P$ is contained in a subarc of $k$ or $c$ that is a wave with respect to $\mathcal E$.

\vskip5pt

\begin{definition} \label{small-large} \rm 
Let $c$ be a simple closed curve on  a surface $\Sigma$ of genus $g \geq 2$. Assume that $c$ is tight 
with respect to some complete decomposing system $\mathcal{E}$ of $\Sigma$. An arc $c'$ on the curve $ c $ will be called {\it small} (with respect to $\mathcal E$) if 
$$| c' |_{\mathcal{E}}\,\, < \,\,   \frac 3 {12g-11} | c |_{\mathcal{E}}\, .$$
If $c'$ is not small, it will be called {\it large} (with respect to $\mathcal E$).

\end{definition}

We now use Definition \ref{tighttts} and assume that $\tau$ is a maximal train track on the surface 
$\Sigma$ which is tight with respect to the complete decomposing system $\mathcal E$. We denote 
by  $| \tau |_\mathcal E$ the {\it $\mathcal E$-length} of $\tau$, by which we mean the total 
$\mathcal E$-length of any set of arcs $\alpha_i$ such that every regular $I$-fiber is met  by only 
one of the $\alpha_i$, and precisely once.

If $\tau'$ is a train track derived from $\tau$, then the length  $| \tau' |_{\mathcal{E}}$  is precisely given by $| \tau |_{\mathcal{E}}$  plus the sum of the lengths $| \sigma_i |_{\mathcal{E}} $ of all of the unzipping paths  $\sigma_i$ used to derive $\tau'$ from $\tau$ (compare subsection \ref{unzippingplusderived}).  
Below we will always use the convention that any unzipping path used to derive $\tau'$ from $\tau$ is oriented from the cusp point of $\tau$ towards the cusp point of $\tau'$.

Since $\tau$ is maximal, each complementary component is a triangle, so that by Euler characteristic reasons there must be precisely $4g-4$ such triangles. Each triangle gives rise to precisely $3$ unzipping paths $\sigma_i$, so that altogether we have $12g -12$  unzipping paths $\sigma_i$.  

\vskip7pt

\begin{remark}\label{new-six-two} \rm 

\noindent (a) By the definition of the deriving process, every $\sigma_i$ must cover $\tau$, 
so that one has: $$|\sigma_i|_\mathcal E \geq |\tau|_\mathcal E$$

\vskip5pt

\noindent
(b) As a consequence, we obtain:
$$| \tau' |_{\mathcal{E}} \, \,= \, \,|\tau|_\mathcal E + \sum_{i = 1, \ldots, 12g - 11}|\sigma_i|_\mathcal E \, \, \geq \, \, (12g - 11) | \tau |_{\mathcal{E}} $$

\end{remark}

\begin{proposition} \label{nonsmall} 

Let $\mathcal E$ be a complete decomposing system of $\Sigma$, let 
$\tau \supset \tau' \supset \tau'' \supset \tau'''$ be a tower of derived train tracks that are tight with respect to $\mathcal E$, and let $c$ be a simple closed curve that is carried by $\tau'''$. Let $c'$ be an arc on $c$ that is large with respect to $\mathcal E$. Then one has:

\begin{enumerate}

\item[(1)]   The arc $c'$ contains a subarc that is parallel with respect to $\mathcal E$ 
to one of the unzipping paths $\sigma_i$ used to derive $\tau'$ from $\tau$.

\item[(2)] The arc $c'$ covers $\tau$.

\item[(3)] Let $\mathcal{D} = \{D_1,\dots, D_{3g-3}\}$ be a second complete decomposing system on 
$\Sigma$ which is tight with respect to $\mathcal E$, and assume that for some intersection point 
$P \in D_i \cap c$  the intersection arc on $c$ is large. Then  $\mathcal D$ can be carried by $\tau$. 

\end{enumerate}

\end{proposition}

\vskip5pt

\begin{proof}
(1) We first observe that, since $c$ is carried by $\tau'''$, it follows from Lemma \ref{carriesimpliesunzip} that $c$ covers $\tau''$.

Since $c'$ is large, we can decompose $c'$ as  a concatenation  $c' = c_1 \circ c_2 \circ c_3$ 
where for each $k \in \{1, 2, 3\}$,  $c_k$ has length 
$|c_k|_\mathcal E \geq \frac 1 {r + 1} | c |_{\mathcal{E}} $. 
Since $c$ covers $\tau''$, we have $| c |_\mathcal{E} \geq | \tau'' |_\mathcal{E}$. Hence 
Remark \ref{new-six-two} (b) yields for any $k \in \{1, 2, 3\}$:

$$|c_k|_\mathcal E \geq \frac 1 {12g-11} | c |_{\mathcal{E}} \geq \frac 1 {12g-11} | \tau'' |_{\mathcal{E}} >  | \tau' |_{\mathcal{E}} $$

In particular, $c_2$ must intersect at least one singular fiber $I_0$ of $\tau'$, say in a point $Q$, and there is at least one of the unzipping paths $\sigma_i$ which has its terminal point on $I_0$. Let $c''$ be the arc on $c$ that starts at $Q$ and runs in the same direction as $\bar \sigma_i$ (= $\sigma_i$ with reverted orientation), and has the same length as $\sigma_i$. From Remark \ref{new-six-two} (a) and 
the above inequality we obtain for $k = 1$ or $k = 3$:
$$ |c''|_\mathcal E = |\sigma_i|_\mathcal E \leq |\tau'|_\mathcal E \leq |c_k|_\mathcal E  $$
Hence $c''$ is a subpath of $c'$. 
 
We now ask whether $c''$ runs parallel to $\sigma_i$ on $\tau'$  (or on $\tau$). The only way in which this can fail to happen is if at some singular fiber $I_1$ of $\tau'$ the two paths branch off each other on $\tau'$. Let $\sigma_j$ be the unzipping path with terminal point at $I_1$. Notice that the two branches of $\tau'$ on either side of $\sigma_j$ run still parallel on $\tau$. Thus either $c''$ and $\sigma_i$ run parallel on $\tau$ throughout all of $\sigma_i$, or else $c''$ runs parallel to all of $\sigma_j$ (or parallel to all of some other $\sigma_k$, which we then rename $\sigma_j$ for the rest of the proof), before branching off $\sigma_i$ on $\tau$. Thus $c''$ runs parallel on $\tau$ for the entire length of either $\sigma_i$, or for the entire length of $\sigma_j$.  Since $\tau$ is tight with respect to $\mathcal E$, the same assertion is true with ``parallel on $\tau$'' replaced by ``parallel with respect to $\mathcal E$''. Hence  assertion (1) of the lemma is proved.

\vskip5pt

\noindent 
(2) This is a direct consequence of the proof given above for assertion (1) since, by definition, 
any of the unzipping paths $\sigma_i$ (and hence also any path parallel to $\sigma_i$ on $\tau$) 
covers $\tau$.

\vskip5pt

\noindent (3) The intersection arc $\alpha'$ at $P$ is an arc on $c$, and since it is large, it follows from assertion (2) that it covers $\tau$. The arc $\alpha$ on $D_i$ which is parallel to $\alpha'$ (compare the paragraph before Definition \ref{small-large}) can be isotoped close to $\alpha'$ and thus to an arc which is parallel to $\alpha'$ on $\tau$, so that it also covers $\tau$. This isotopy preserves the property that $\mathcal D$ is tight with respect to $\mathcal E$, as by definition (see the beginning of this section) it is an isotopy of the pair $(\Sigma, \mathcal{E})$. We can thus apply Proposition \ref{carriedcurves}
and obtain that $\mathcal D$ is carried by $\tau$.

\end{proof}

\begin{proposition} \label{newdichotomy}
Let $\mathcal E \subset \Sigma$ be a complete decomposing system and let $c$ and $D$ be 
simple closed essential curves on $\Sigma$ which are tight with respect to  $\mathcal{E}$. We also assume that for any intersection point  $P \in D \cap c$ the intersection length $|P|_{\mathcal{E}}$ is  small. 

Let $\delta_c^m(D)$ denote the curve obtained from $D$ by an $m$-fold Dehn twist at $c$ and 
subsequent  tightening with respect to $\mathcal E$. Then for any non-trivial twist exponent 
$m \in \mathbb Z \smallsetminus \{ 0 \}$ and any intersection point $S \in \delta_c^m(D) \cap c$ 
the intersection length $|S|_\mathcal E$ is large.

\end{proposition}

\begin{proof}  At every intersection point $P \in D \cap c$ we perform the $m$-fold Dehn twist  at $P$
in a two-step procedure as follows:

\vskip5pt

\noindent \underline{Step 1:} Choose a small embedded annulus neighborhood $A_P$ of $c$ in 
$\Sigma$. The point $P$ is contained  in  the arc  $\beta = A_P \cap D$.  After a suitable isotopy of $D$ we can assume that $\beta$ is entirely contained in one  of the pair-of-pants complementary to $\mathcal E$. Denote the points in $\partial \beta$ by $P^{in}$ 
and $P^{tr}$ (``initial'' and ``terminal'' points).  Remove the arc $\beta$ from $ D$ and instead insert  
an arc $\eta \subset A_P$ with $\partial \eta = \{P^{in}, P^{tr}\}$, where $\eta$ winds around $m$ times around the core curve of $A_P$ in the direction determined according to whether $m > 0$ or $m < 0$.
After a suitable isotopy we can assume that the arc $\eta$ meets $c$ only in a single point $S$.

\vskip5pt

\noindent \underline{Step 2:} Now perform an isotopy of the new curve 
$\hat{D}  = (D \smallsetminus \beta) \cup \eta$, which tightens it with respect to $\mathcal{E}$, in 
order to obtain the curve $\delta^m_c(D)$. This is done by isotoping off $\mathcal{E}$ two pairs of parallel arcs: The arcs in the first pair are  concatenated at $P^{in}$  and the arcs in the second 
pair are concatenated at $P^{tr}$. In each pair  one of the arcs lies on $\eta$, and the other on 
$D \smallsetminus \beta$.  The tightening isotopy will move the points $P^{in}$, $P^{tr}$ in opposite directions along  paths determined by each pair of parallel arcs.

\vskip5pt

There are now several cases to be considered, and in order to do so in a precise way, we 
introduce  the following notation for the arcs in $D$ that will be cancelled as described in the 
above step 2: 

Consider the intersection arc $\alpha$ of $P$ on $ D$. Since $c$ and $D$ are tight with respect 
to $\mathcal E$, the point $P$ lies on the arc $\hat{\alpha}$ which is the prolongation of $\alpha$ 
into the pairs-of-pants adjacent to the curves of $\mathcal{E}$ containing $\partial \alpha$. (Recall 
that $\partial \alpha \subset \mathcal{E}.)$ Denote by $\hat{D}^{in}$ and  $\hat{D}^{tr}$ the two connected components of $\alpha - \beta$.  Recall that all of $\beta$ and hence $P, P^{in}$ and 
$P^{tr}$ are contained in the interior of one of  the pair-of-pants complementary to $\mathcal E$. 
In case where $\beta$ is contained in  $\hat \alpha \smallsetminus \alpha$, 
then one of the two, 
$\hat{D}^{in}$ or $\hat{D}^{tr}$, is empty. Similarly, denote by $\eta^{in}$ and $\eta^{tr}$ the 
maximal initial and terminal subarcs of $\eta$ 
which run parallel to $\hat{D}^{in}$ and $\hat{D}^{tr}$ respectively. Note that the arcs  $\eta^{in}$ 
and $\hat{D}^{in}$ are concatenated at $P^{in}$ and the arcs  $\eta^{tr}$ and $\hat{D}^{tr}$ are concatenated at $P^{tr}$. Note also  that both pairs of  concatenated arcs can be cancelled  by an isotopy that moves the points $P^{in}$ and $P^{tr}$ in opposite directions along $\hat{\alpha}$. Furthermore, note the following crucial fact: 

\vskip5pt

\noindent
($\ast$) By the definition of the intersection arc $\alpha$ the cancelling isotopy 
can not be extended further, as the extension of the above arcs, along $\eta$
and $D$ respectively, must leave each of the two pairs-of-pants adjacent to $\alpha$
through distinct boundary curves.

\vskip5pt

If there is only a single intersection point $P$ of $D$ with $c$, then the above paragraph gives a precise description of an $m$-fold Dehn-twist of $D$ along  $c$. However, if there are several 
such intersection points, one has to be  much more careful: The difficulty comes from the fact that for distinct  intersection points the above subarcs may overlap, so that it becomes  impossible to perform 
the above tightening isotopies at all intersection  points at the same time. One also needs an argument to show that after  performing at adjacent intersection points two such isotopies which end  in the same pair-of-pants complementary to $\mathcal{E}$, there is no possibility  to continue the tightening process further. 

In order to analyse the cases systematically, we first observe that for a single intersection point $P$ the subarcs $\eta^{in}$ and $\eta^{tr}$ of $\eta$ can not overlap (on $\eta$): This is due to the hypothesis that the intersection length at $P$ is small, since both $\eta^{in}$ and $\eta^{tr}$ must run parallel to a subarc of the intersection arc at $P$, while the length of $\eta$ is a multiple of the length of $c$.

Hence for each point $P$  we can choose a point $\bar P$ on $\mathcal{E} \cap \eta$ which is 
not contained in either $\eta^{in}$ or $\eta^{tr}$. We will now show that,  if  $P$ and $P'$ are 
subsequent intersection points on $D$,  then the segment  $[\bar P, \bar P'] \subset \hat D$ between 
the corresponding points $\bar P$ and  $\bar P'$ on $\hat D$  becomes tight after cancelling the segments $\eta^{in}$, $\eta^{tr}$,  $\hat D^{in}$ and $\hat D^{tr}$, or subsegments of the latter, but that no  further cancelation of  intersection points with $\mathcal E$ ever occurs in the process of
tightening the segment $[\bar P, \bar P']$.

Note that, if the intersection length at $P$ is bigger than $0$, one can isotope the intersection point 
$P$, and with it the arc $\beta$, along the intersection arc of $P$ on $c$. Such an isotopy gives rise 
to a ``trade-off'' between the lengths  of $\eta^{in}$  and $\hat D^{in}$ on one hand, and  $\eta^{tr}$ 
and $\hat D^{tr}$ on the other.  Of course, the total number of possible cancellations is not 
affected by such a trade-off move. Furthermore, the  corresponding point $\bar P$ introduced  
above can be kept fixed while performing a trade-off move at $P$.

As before, let  $P$ and $P'$ two intersection points that are adjacent on $D$, and let  $\hat d$ 
denote the segment of $D$ which lies between the subarcs $\beta$ and $\beta'$ corresponding 
as above to $P$ and $P'$. Since $P$ and $P'$ are adjacent, $\hat d$ does not contain another intersection point of $D$ with $c$. Thus $\hat d$ can also be viewed as subarc of  $\hat D$: It is precisely the segment between the subarcs $\eta$ and $\eta'$.
We distinguish three cases:

\vskip5pt
\noindent
(1) The intersection arcs of $P$ and $P'$ on $D$ do not overlap along $\hat d$, and furthermore 
they  are separated on $\hat d$ by at least one intersection point, say $R$, of $D \cap \mathcal E$.
Thus in this case, the corresponding cancellation arcs on $\hat D$, say $\hat D^{tr}$ and  
$\hat {D'}^{in}$,  cannot overlap on $\hat d$,  as they are separated by $R$.  Hence, by the fact ($\ast$) stated above, after canceling $\hat D^{tr}$ against $\eta^{tr}$ and $\hat {D'}^{in}$ against $\eta'^{in}$ the resulting segment between  $\bar P$ and   $\bar P'$  will be tight (see Figure 1).
 
 \vskip5pt

\begin{figure}
{\epsfxsize = 5.0 in \centerline{\epsfbox{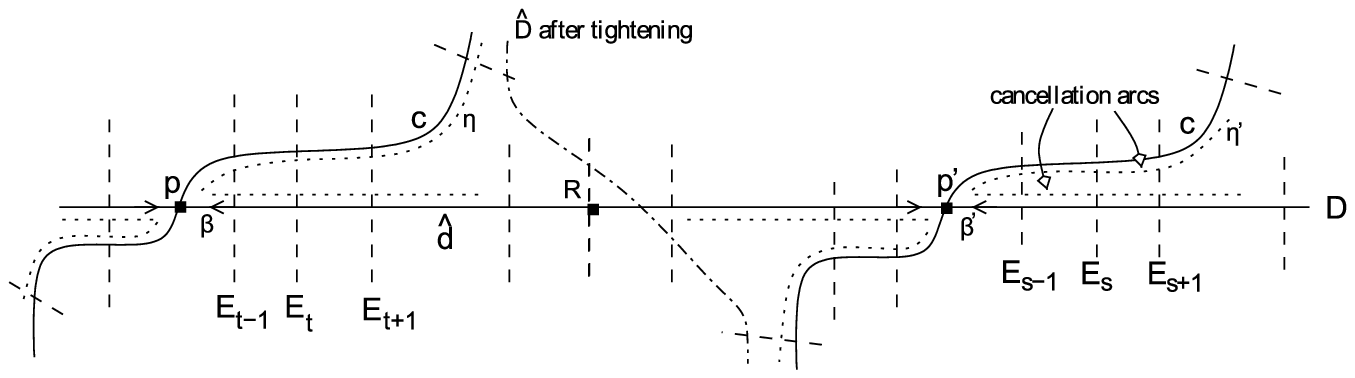}}}
\caption{The intersection arcs do not overlap over $\hat{d}$.}

\end{figure}

 \vskip5pt
 
\noindent
(2)  The intersection arcs of $P$ and $P'$ on $\hat D$ do overlap along $\hat d$. Hence we can 
perform a trade-off move as defined above, so that after this isotopy both, $P$ and $P'$, together 
with  the small corresponding arcs $\beta$ and $\beta'$, all come to lie in the same complementary 
pair-of-pants $\mathcal P$, and $\hat{d}$ becomes a very small arc contained in $\mathcal P$ (and 
thus disjoint from $\mathcal E$). Note that this can be done only if the annulus neighborhood $A$ 
of $c$ was chosen sufficiently thin. Now observe that the directions of the two inserted arcs $\eta$ 
and $\eta'$ coincide, as the ``twisting direction'' of a Dehn twist is well defined and independent of the local orientation of the curve $c\,$: This is a well known fact for Dehn twists.

As a consequence, in this second case no cancellation at all is possible on the segment of 
$\hat D$ between $\bar P$ and $\bar P'$: This arc is parallel to  a concatenation of $\hat d$ 
with two arcs that each winds around $c$ at most $m$ times, 
and thus is already tight as is     (see Figures 2 and 3).

 \vskip0pt
 
\begin{figure}
{\epsfxsize = 5.0 in \centerline{\epsfbox{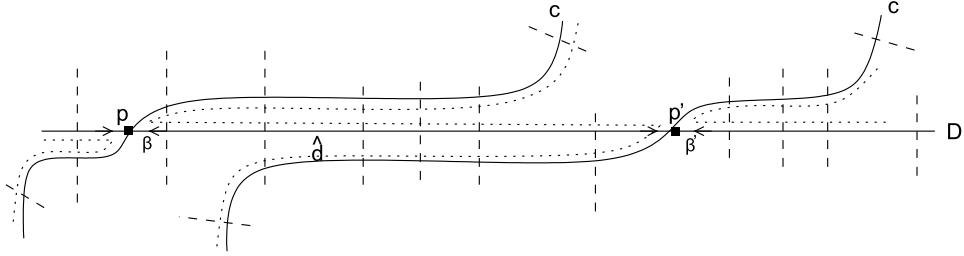}}}
\caption{The intersection arcs do overlap over $\hat{d}$.}

\end{figure}
 
 \vskip0pt
 
 \begin{figure}
{\epsfxsize = 5.0 in \centerline{\epsfbox{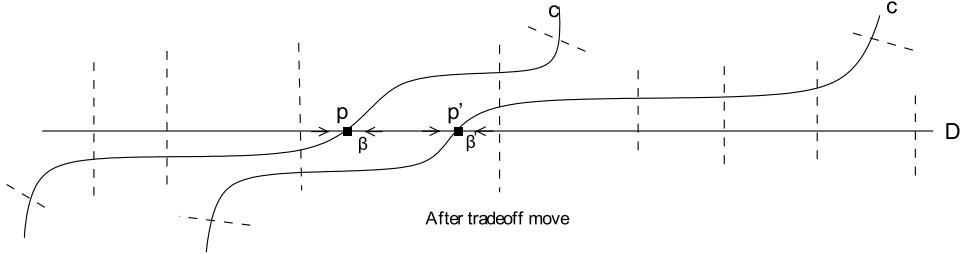}}}
\caption{Itersection arcs after the ``trade off''.}

\end{figure}

\vskip2pt

\noindent (3) The intersection arcs of $P$ and $P'$ on $\hat D$ do not overlap along $\hat d$, but they have endpoints $Q$ and $Q'$ which are contained in curves $E_i, E_j \in \mathcal{E}$, and 
$Q, Q'$ are adjacent points of $\hat D \cap \mathcal{E}$ on $\hat d$:

This implies that there is  a pair-of-pants $\mathcal P$ complementary to $\mathcal E$ which contains both, $Q$ and $Q'$, but on distinct boundary curves of $\mathcal P$. In this case we can again perform a trade-off move, so that $P$ and $P'$ move along $\hat d$ beyond $Q$ and $Q'$ respectively, and both come to lie inside of $\mathcal P$.

But then the same argument, as in case (2) above, applies: The  arcs $\beta$ and $\beta'$ can 
also be assumed to lie in $\mathcal P$, the segment $\hat d$ degenerates to a small arc entirely contained in $\mathcal P$, and the segment on $\hat D$ between $\bar P$ and $\bar P'$ 
is tight as is, without any cancellations at either  $\hat P^{tr}$ or $\hat P'^{in}$ (see Figure 4).

 \vskip5pt
 
 \begin{figure}
{\epsfxsize = 5.0 in \centerline{\epsfbox{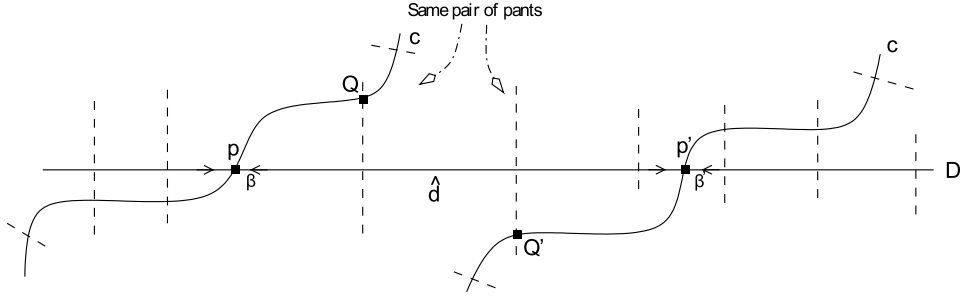}}}
\caption{End points Q and Q' on boundary curves of the same pair of pants .}

\end{figure}

\vskip10pt

Hence in each of the three possible cases the tightening process needed to obtain  $\delta_c^m(D)$ 
from $\hat D$ is complete if at any one of the intersection points  $P \in D \cap \mathcal E$ an isotopy 
is performed which cancels $\eta^{in}$ against  $\hat D^{in}$ and $\eta^{tr}$ against $\hat D^{tr}$, or subarcs of those. 

However, by definition  the intersection arc $\alpha$ of $P$ on $D$ is equal to the concatenation 
$\hat D^{in} \circ \beta \circ \hat D^{tr}$, and $\beta$ is disjoint from $\mathcal E$.  Hence one obtains 
$|P|_\mathcal E = |\alpha|_\mathcal E = |\hat D^{in}|_\mathcal E + |\hat D^{tr}|_\mathcal E = |\eta^{in}|_\mathcal E + |\eta^{tr}|_\mathcal E$.

On the other hand, all of the arc $\eta$ runs parallel to an arc on $c$, and one has 
$|\eta|_\mathcal E = m |c|_\mathcal E$. Hence the intersection length of $\delta_c^m(D)$ at $S$ is 
bounded below by $|\eta \smallsetminus (\eta^{in} \cup \eta^{tr})|_\mathcal E = |\eta |_\mathcal E - (|\eta^{in}|_\mathcal E + | \eta^{tr})|_\mathcal E)$.
For $|m| \geq 1$ we obtain:

\vskip0pt

$$|S|_\mathcal{E} \, \, \geq \, \,  |\eta |_\mathcal E - (|\eta^{in}|_\mathcal E + | \eta^{tr})|_\mathcal E) \,\, \geq \, \, |m| \cdot |c|_{\mathcal{E}} -  |P|_{\mathcal{E}} \, \, $$ 

$$> \, \, |c|_{\mathcal{E}} - \frac 3 {12g - 11} | c |_{\mathcal{E}} \,\, 
>\, \, \frac {r-2} {12g - 11} | c |_{\mathcal{E}} \,\, >\, \,  \frac {12g-14} {12g - 11} | c |_{\mathcal{E}} \,\, >\, \, 
\frac {3} {12g - 11} | c |_{\mathcal{E}} $$

\vskip20pt

This shows that the intersection length of $\delta_c^m(D)$ at any of its intersection points with $\mathcal E$ is large.

\end{proof}

\vskip10pt

\begin{corollary}
\label{minimum}
Let $\mathcal E \subset \Sigma$ be a complete decomposing system and let $c$ and $D$ be 
simple closed essential curves on $\Sigma$ which are tight with respect to  $\mathcal{E}$. We also assume that for every intersection point  $P \in D \cap c$ the intersection length $|P|_{\mathcal{E}}$ is  small. 

Let $\delta_c^m(D)$ denote the curve obtained from $D$ by an $m$-fold Dehn twist at $c$ and 
subsequent  tightening with respect to $\mathcal E$. Then   the total number of intersections of 
$\delta_c^m(D)$ with $\mathcal E$ achieves a strict minimum at $m = 0$, where  $m \in \mathbb Z$
is the twist exponent:
$$\# \delta_c^m(D) \cap \mathcal E \,\, > \, \, \# D \cap \mathcal E \qquad {\rm for \, all} \qquad  m \neq 0$$

\end{corollary}

\begin{proof}
This is a direct consequence of the fact, shown in detail in the last proof, that the curves $D$ and $\delta_c^m(D)$ differ only along the intersection arcs of their intersection points with $c$, and for those the statement of Proposition \ref{newdichotomy} gives directly the desired inequality.

\end{proof}

\vskip10pt

We will now go back to the setting of subsection \ref{subsection-Heegaard-splittings}:
Let $M$ be a closed $3$-manifold with Heegaard splitting $M = V  \cup_{\Sigma} W$, and let
$n \in \naturals$.  A curve $c$ on the Heegaard surface  $\Sigma$ is called {\it $n$-optimal} if for 
any $m \in  \integers \smallsetminus \{0\}$ one has
$$ d(V, W_{c}^m) \geq n \, ,$$
where  the handlebody $W_c^m$ is obtained from $W$ by composing the attaching map 
$\partial W \to \Sigma$ with an $m$-fold Dehn  twist along $c$ (compare the discussion at the 
beginning of section \ref{introduction}). The induced Heegaard splitting $(V, W_{c}^m)$ defines 
a 3-manifold $M^m_c = V \cup_\Sigma W_{c}^m$  which can be obtained alternatively by horizontal surgery at $c \subset \Sigma$, see Remark \ref{threesphere}.

\vskip10pt

\begin{definition} \label{gregariousnew} \rm
A simple closed curve $c$ on $\Sigma$ is called {\it $n$-gregarious} with respect to some complete decomposing system $\mathcal E$ on $\Sigma$, for any integer $n \geq 0$, if $c$ is $n$-gregarious with respect to some fat train track $\tau$ on $\Sigma$ with $\mathcal{E}_\tau = \mathcal E$. 

\end{definition}

\vskip10pt

\begin{theorem} \label{surgerythm}
Let $M$ be an oriented $3$-manifold with a Heegaard splitting  $M = V   \cup_{\Sigma} W$, 
and let  $n \in \mathbb{N}$ be an integer that satisfies $n > d(V, W)$.  Consider any complete decomposing systems  $\mathcal D \in  \mathcal{CDS}(V)$ and   $\mathcal E \in \mathcal{CDS}(W)$.
Then any essential simple closed  curve $c$ on $\Sigma$ that is  $(n+3)$-gregarious with respect to both, $\mathcal D$ and  $\mathcal E$,  is $n$-optimal.

\end{theorem}

\vskip5pt

\begin{proof}
(a)
We first use Corollary \ref{minimum} to argue that if for some exponent $m_0 \in \mathbb Z$ the intersection length $|P|_\mathcal E$ is small, at every intersection point $P \in D_i \cup c$ and for 
every $D_i \in \delta_c^{m_0}(\mathcal D)$, then the total number of intersection points satisfies
$$\#(\delta_c^{m}(\mathcal{D}) \cap \mathcal E) \,\, > \,\, \#(\delta_c^{m_0}(\mathcal{D}) \cap \mathcal E)$$
for all exponents $m \neq m_0$. Since the intersection number is invariant under homeomorphisms, we can apply $\delta_c^{-m}$ to obtain 
$$\# (\delta_c^{m}(\mathcal{D}) \cap \mathcal E) \,\, = \,\, \# (\mathcal D \cap \delta_c^{-m}(\mathcal{E}))$$
for all $m \in \mathbb Z$. From the symmetry of the conditions in $\mathcal D$ and $\mathcal E$ we deduce that, if for some exponent $m_1 \in \mathbb Z$ the intersection length $|Q|_\mathcal D$ is small, at every intersection point $Q \in E_j \cup c$ and for every $E_j \in \delta_c^{m_1}(\mathcal E)$, then one must have:
$$m_0 = m_1$$

\vskip5pt

\noindent
(b)
We now consider a situation where for some $m \in \mathbb Z$ there are intersection points  $P \in D_i \cup c$ for some $D_i \in \delta_c^{m}(\mathcal D)$ and $Q \in E_j \cup c$ for some $E_j \in \delta_c^{m}(\mathcal E)$ such that both, $|P|_\mathcal E$ and $|Q|_\mathcal D$, are large. In this case we can apply Proposition \ref{nonsmall} (c) to deduce that $\mathcal D$ and $\mathcal E$ are $n$-gregarious with respect to each other, and hence in particular they don't have waves with respect to each other.  Thus we can apply Theorem \ref{HSdistance} to deduce that $d(V, W_c^m) \geq n$. As a consequence, it follows from the hypothesis $n > d(V, W)$ that $m \neq 0$.

\vskip5pt

\noindent
(c)
It follows from the above part (a) that there is at most one value $m_0 \in \mathbb Z$ for which no point $P$ as in the hypothesis of part (b) exists, and at most one $m_1 \in \mathbb Z$ without an analogous point $Q$. Furthermore, from (a) we know that if there exist both, $m_0$ and $m_1$, then one has $m_0 = m_1$.  In any case, there is at most one exceptional value $m' = m_0, m' = m_1$ or $m' = m_0 = m_1$, 
such that the hypotheses of part (b) are given for all integers $m \neq m'$. Thus we conclude from part (b) that $d(V, W_c^m) \geq n$ for all $m \in \mathbb{Z} \smallsetminus \{m'\}$, and that hence that $m' = 0$. This  proves that the curve $c$ is $n$-optimal.

\end{proof}

\vskip10pt

\begin{theorem} \label{twistgenericity}
Let $M$ be a closed $3$-manifold with Heegaard splitting $M = V  \cup_{\Sigma} W$. For any  
$n \geq 1$ set of $\mathcal{C}^{M}_{n}(\Sigma)$ of $n$-optimal curves is  generic in the set 
$\mathcal{C}^{0}(\Sigma)$ of all essential simple closed curves on $\Sigma$.

\end{theorem}

\vskip10pt

\begin{proof}
We first observe that for any integer $n' \geq n$ one has $\mathcal{C}^{M}_{n'}(\Sigma) \subset \mathcal{C}^{M}_{n}(\Sigma)$, by the definition of $n$-optimal. Thus, by Remark \ref{reformulation} (b), it suffices to prove the statement for $n > d(V, W)$.

We apply Lemma \ref{quotehempel} to find $\mathcal D$ and $\mathcal E$ as in Theorem \ref{surgerythm}, and deduce from the latter that the set  $\mathcal{C}^{M}_{n}(\Sigma)$ contains the intersection of the set  ${\bf G^{n+3}}(V)$ and the set   ${\bf G^{n+3}}(W)$. According to  Proposition 
\ref{morelongdistance} both of these sets are generic in $\mathcal{C}^{0}(\Sigma)$. Hence part (a) of 
Remark \ref{reformulation} shows  that the intersection is generic, and part (b) implies that the set 
$\mathcal{C}^{M}_{n}(\Sigma)$ is generic in $\mathcal{C}^{0}(\Sigma)$.

\end{proof}

\end{document}